\documentclass[a4paper,12pt]{amsart}
\usepackage{fullpage}
\usepackage{amsmath,amssymb,mathtools,booktabs,tikz}
\usetikzlibrary{arrows.meta,calc,decorations.markings,positioning}
\usepackage[T1]{fontenc}
\usepackage{lmodern}
\usepackage[expansion=false]{microtype}
\usepackage{needspace}
\usepackage{tabularx,array}
\usepackage[hypertexnames=false]{hyperref}
\usepackage{doi}
\numberwithin{equation}{section}
\usepackage{float}

\usepackage{graphicx}
\definecolor{potfigink}{HTML}{233B48}
\definecolor{potfigblue}{HTML}{167D9A}
\definecolor{potfigrust}{HTML}{B45635}
\definecolor{potfigpurple}{HTML}{80649C}
\definecolor{potfigbluewash}{HTML}{E8F3F5}
\definecolor{potfigrustwash}{HTML}{FBEEE7}
\definecolor{potfigcollar}{HTML}{D6E9D5}
\definecolor{potfigmuted}{HTML}{64737B}
\definecolor{potfiggrid}{HTML}{DCE2E5}
\hypersetup{colorlinks=true,linkcolor=blue,citecolor=blue,urlcolor=blue,
  pdftitle={Equilibrium points of logarithmic potentials},
  pdfauthor={Alexandre Eremenko and Teng Zhang}}

\newtheorem{thm}{Theorem}[section]
\newtheorem{lem}[thm]{Lemma}
\newtheorem{prop}[thm]{Proposition}
\newtheorem{cor}[thm]{Corollary}

\newtheorem*{thmA}{Theorem A}
\newtheorem*{thmB}{Theorem B}
\newtheorem*{thmC}{Theorem C}
\theoremstyle{remark}

\makeatletter
\newcommand{\namedtheoremlabel}[2]{%
  \phantomsection
  \def\@currentlabel{#2}%
  \def\@currentlabelname{Theorem #2}%
  \label{#1}%
}
\makeatother
\newcommand{\C}{\mathbb C}
\newcommand{\bC}{\overline{\mathbb C}}
\newcommand{\R}{\mathbb R}
\newcommand{\Z}{\mathbb Z}

\newcommand{\HH}{\mathbb H}
\newcommand{\dd}{\,d}
\newcommand{\Rea}{\operatorname{Re}}
\newcommand{\Ima}{\operatorname{Im}}
\newcommand{\sgn}{\operatorname{sgn}}
\newcommand{\Log}{\operatorname{Log}}
\newcommand{\smin}{s_{\min}}

\DeclareMathOperator{\dist}{dist}
\title[Equilibrium points of logarithmic potentials]
{Equilibrium points of logarithmic potentials}
\author[A.~Eremenko]{Alexandre Eremenko}

\address{Mathematics Department, Purdue University,
	West Lafayette, IN 47907, USA}
\email{eremenko@purdue.edu}

\author[T.~Zhang]{Teng Zhang}
\address{School of Mathematics and Statistics, Xi'an Jiaotong University,
Xi'an 710049, P. R. China}
\email{teng.zhang@stu.xjtu.edu.cn}
\subjclass[2020]{Primary 31A15; Secondary 30C62, 30D35}
\keywords{Logarithmic potential, critical point, Cauchy transform,
quasiconformal mapping, analytic capacity}
\begin{document}
\begin{abstract}
We construct logarithmic potentials with positive masses at a sequence
of points tending to infinity and with no critical points. These
potentials can have any prescribed order $\rho\in[1/2,1)$; for
$\rho=1/2$, we obtain examples of both normal and maximal type.
For $1/2<\rho<1$, the masses at conjugate poles are prescribed by a
power law. We determine the pole asymptotics and the angular growth
of these potentials. For Cauchy sums with absolutely summable complex
residues and nonzero total residue, we prove the existence of infinitely
many zeros and compare their counting function with that of the poles.
\end{abstract}
\maketitle
\raggedbottom

\section{Introduction}

Consider positive point masses placed at a sequence of points in the
plane tending to infinity. Must their forces balance somewhere?
Clunie, Eremenko and Rossi~\cite{CER93} asked this question for
logarithmic and Newtonian potentials, extending a problem of Rubel
collected by Brannan and Hayman~\cite[Problem~7.78]{BH89}.

In the complex plane, let $(\zeta_j)_{j\ge1}$ be an infinite sequence
of distinct nonzero complex numbers, and assign a positive mass $a_j$
to each site $\zeta_j$. Assume that
\begin{equation}\label{eq:sum}
 |\zeta_j|\to\infty\quad(j\to\infty),\qquad
 \sum_{j=1}^{\infty}\frac{a_j}{|\zeta_j|}<\infty.
\end{equation}
The associated logarithmic potential, normalized by $u(0)=0$, is
\begin{equation}\label{eq:cauchy-potential}
 u(z)=\sum_{j=1}^{\infty}a_j
      \log\left|1-\frac{z}{\zeta_j}\right|,
\end{equation}
and its conjugate gradient is
\begin{equation}\label{mero}
 f(z)=u_x(z)-iu_y(z)
     =\sum_{j=1}^{\infty}\frac{a_j}{z-\zeta_j}.
\end{equation}
By \eqref{eq:sum}, both series converge absolutely and locally
uniformly away from the sites. Thus $f$ is meromorphic, with simple
poles at $\zeta_j$ and positive residues $a_j$. We call $f$ a
\emph{positive Cauchy sum}. The critical points of $u$ are precisely
the zeros of $f$.

The potential $u$ is subharmonic in the plane. Its order is
\[
 \rho[u]=\limsup_{r\to\infty}\frac{\log B(r,u)}{\log r},
 \qquad B(r,u)=\max_{|z|\le r}u(z).
\]
Condition~\eqref{eq:sum} implies $B(r,u)=o(r)$: indeed,
$u^+(z)\le\sum_j a_j\log(1+|z|/|\zeta_j|)$, and division by
$|z|$, followed by dominated convergence, proves the assertion.
Consequently, $\rho[u]\le1$, and $u$ is of minimal type if
$\rho[u]=1$. We also write
\[
 M(r)=\sum_{|\zeta_j|\le r}a_j
\]
for the total mass in the disk $|z|\le r$.

We use the standard notation of Nevanlinna theory as in~\cite{GO08}.
The order and lower order of a meromorphic function $f$ are defined by
the limsup and liminf, respectively, of $\log T(r,f)/\log r$.
 For a Cauchy sum satisfying
\eqref{eq:sum}, Lemma~\ref{lem:proximity} shows that its order is
determined by the pole sequence.

The existence and distribution of critical points of such potentials
have been studied in~\cite{CER93,ELR94,LR04,Lan07,She24}. We quote two results.

\begin{thmA}[{\cite[Theorem~2.10]{CER93}}]
\namedtheoremlabel{thm:A}{A}
Let $u$ be the positive logarithmic potential in
\eqref{eq:cauchy-potential}, and assume that \eqref{eq:sum} holds.
If $B(r,u)=o(\sqrt r)$, then $u$ has infinitely many critical points.
\end{thmA}

Theorem~\ref{thm:A} is stated under the condition $M(r)=o(\sqrt r)$.
The formulation above follows from Jensen's formula, which gives
$M(r)\log2\le B(2r,u)$.

\begin{thmB}[{\cite[Corollary~1.4]{ELR94}}]
\namedtheoremlabel{thm:B}{B}
Let $u$ be as in \eqref{eq:cauchy-potential}, and assume that
\eqref{eq:sum} holds. If $\inf_{j\ge1}a_j>0$, then $u$ has infinitely
many critical points.
\end{thmB}

Clunie, Eremenko and Rossi~\cite[Conjecture~2.7]{CER93} conjectured
that every positive Cauchy sum satisfying \eqref{eq:sum} has a zero.
In view of results such as Theorem~\ref{thm:B}, Langley and
Rossi~\cite[Conjecture~1.1]{LR04} stated the stronger conjecture
that every such sum has infinitely many zeros. Our examples disprove
both conjectures and show that the growth condition in
Theorem~\ref{thm:A} is sharp.

We begin with potentials of order $\rho\in(1/2,1)$.
Put $\beta=2\rho-1\in(0,1)$ and prescribe the masses
\begin{equation}\label{eq:mass-profile}
 b_m=m^{\beta-1},\qquad m\ge1.
\end{equation}
Each mass $b_m$ will be placed at both points of a conjugate pair
$z_m,\overline{z_m}$, whose locations are to be determined. Put
\begin{equation}\label{eq:main-exponents}
 \rho=\frac{1+\beta}{2},\qquad
 p=\frac{\beta}{\rho}=2-\frac1\rho,\qquad
 \theta=\frac{\pi p}{2}=\pi\left(1-\frac1{2\rho}\right).
\end{equation}

\begin{thm}\label{thm:main}
For the sequence $(b_m)_{m\ge1}$ in \eqref{eq:mass-profile},
there exist a sequence $(z_m)_{m\ge1}$ of distinct points in the
upper half-plane and a constant $\kappa>0$ such that
\[
 z_m\sim\kappa m^p e^{i\theta}\qquad(m\to\infty),
\]
and the Cauchy sum
\begin{equation}\label{eq:main-series}
 f(z)=\sum_{m\ge1}b_m
      \left(\frac1{z-z_m}+\frac1{z-\overline{z_m}}\right)
\end{equation}
converges absolutely away from its poles and has no zeros.
\end{thm}

The pole asymptotics give $b_m/|z_m|=O(m^{\beta-1-p})$.
Since $p>\beta$, we have $\sum_m b_m/|z_m|<\infty$.
Thus \eqref{eq:main-series} is the conjugate gradient of a
logarithmic potential of the form \eqref{eq:cauchy-potential}.

To describe the angular growth, define
\begin{equation}\label{eq:angular-profile}
 h_\rho(z)=|z|^\rho\max\{0,\cos(\rho\arg(-z))\},
 \qquad h_\rho(0)=0,
\end{equation}
where $\arg(-z)\in[-\pi,\pi]$ for $z\ne0$. The definition is
independent of the choice of argument on the cut. Moreover,
$h_\rho(z)=0$ precisely when $z=0$ or $|\arg z|\le\theta$.

\begin{cor}\label{cor:growth}
Let $u$ be the potential associated with the Cauchy sum in
Theorem~\ref{thm:main}, and put $D_\rho=2\kappa^{-\rho}/\beta$.
Then $M(r)\sim D_\rho r^\rho$ and
\begin{equation}\label{eq:angular-limit}
 r^{-\rho}u^+(re^{it})\longrightarrow\pi D_\rho h_\rho(e^{it})
 \quad\text{uniformly for }t\in[-\pi,\pi].
\end{equation}
The convergence
\begin{equation}\label{crr}
 r^{-\rho}u(re^{it})\longrightarrow\pi D_\rho h_\rho(e^{it})
\end{equation}
also holds in $L^1[-\pi,\pi]$. Furthermore,
$B(r,u)\sim\pi D_\rho r^\rho$, while the order and lower order
of $f$ are both $1/p=\rho/(2\rho-1)$.
\end{cor}

The limits in \eqref{eq:angular-limit} and \eqref{crr} describe
completely regular growth in the subharmonic formulation of
Levin--Pfl\"uger theory; see~\cite{Lev80}. More explicitly,
$r^{-\rho}u(rz)\to\pi D_\rho h_\rho(z)$ in $L^1_{\mathrm{loc}}(\C)$.
This follows by polar integration from \eqref{crr} and the bound
$B(r,u)=O(r^\rho)$, using the nonnegative circular mean of $u$.
The angular indicator is $\pi D_\rho h_\rho(e^{it})$; we call
$h_\rho(e^{it})$ the normalized angular profile.

The next theorem treats the endpoint $\rho=1/2$.
At this order, normal type means
$0<\limsup_{r\to\infty}B(r,u)/\sqrt r<\infty$; maximal type
means that this limsup is infinite while the order remains $1/2$.
Put $L=2\pi$ and choose $h$ and $S$ according to the table
\begin{equation}\label{eq:endpoint-models}
 \begin{array}{c|c|c}
  \text{type}&h(t)&S(r)\\ \hline
  \text{normal}&\log t&\sqrt r\\
  \text{maximal}&\exp\sqrt{\log t}&\sqrt{r\log r/2}
 \end{array}.
\end{equation}
Here $h$ is defined for $t\ge L$, and $S(r)$ is considered for
sufficiently large $r$. For each choice, define
\begin{equation}\label{eq:endpoint-masses}
 b_1=\frac{h(L)}{L},\qquad
 b_m=\frac{h(Lm)-h(L(m-1))}{L}\quad(m\ge2).
\end{equation}

\begin{thm}\label{thm:endpoint}
For each choice in \eqref{eq:endpoint-models}, let $(b_m)_{m\ge1}$
be defined by \eqref{eq:endpoint-masses}. There exists a sequence
$(z_m)_{m\ge1}$ of distinct points in the upper half-plane, with
$|z_m|\to\infty$ and $\sum_{m\ge1}b_m/|z_m|<\infty$, such that
the logarithmic potential
\[
 u(z)=\sum_{m=1}^{\infty}b_m
       \left(\log\left|1-\frac{z}{z_m}\right|
       +\log\left|1-\frac{z}{\overline{z_m}}\right|\right)
\]
has no critical points outside its poles. The upper half-plane poles
can be chosen so that, as $m\to\infty$,
\begin{equation}\label{eq:endpoint-poles}
 |z_m|\sim
 \begin{cases}
  (\log m)^2,&\text{normal case},\\
  \exp(2\sqrt{\log m})/\sqrt{\log m},&\text{maximal case}.
 \end{cases}
\end{equation}
Moreover, as $r\to\infty$,
\begin{equation}\label{eq:endpoint-angular}
 \frac{u^+(re^{it})}{S(r)}\longrightarrow|\sin(t/2)|
 \quad\text{uniformly for }t\in[-\pi,\pi].
\end{equation}
The functions $u(re^{it})/S(r)$ also converge to $|\sin(t/2)|$
in $L^1[-\pi,\pi]$. The total mass $M(r)$ and the maximum $B(r,u)$
satisfy
\begin{equation}\label{eq:endpoint-means}
 M(r)\sim S(r)/\pi,\qquad B(r,u)\sim S(r).
\end{equation}
\end{thm}

The normal-type example has $B(r,u)\sim\sqrt r$ and no critical
points, so the condition in Theorem~\ref{thm:A} cannot be weakened
to $B(r,u)=O(\sqrt r)$.

The examples in Theorems~\ref{thm:main} and~\ref{thm:endpoint}
have infinite total mass. For finite total mass, we recall the
classical theorem recorded by Gol'dberg and
Ostrovskii~\cite[Chapter~V, Theorem~6.2]{GO08}. The formulation below,
including finite lower order, is stated in~\cite[Theorem~1.1]{LR04}.

\begin{thmC}[{\cite[Theorem~1.1]{LR04}}]
\namedtheoremlabel{thm:C}{C}
Let $(\zeta_j)_{j\ge1}$ be distinct complex numbers tending to
infinity, and let $(a_j)_{j\ge1}$ be nonzero complex numbers with
\begin{equation}\label{go}
 \sum_{j=1}^{\infty}|a_j|<\infty,\qquad
 A=\sum_{j=1}^{\infty}a_j\ne0.
\end{equation}
If $f(z)=\sum_{j\ge1}a_j/(z-\zeta_j)$ has finite lower order,
then $\delta(a,f)=0$ for every $a\in\bC$.
In particular, $f$ has infinitely many zeros.
\end{thmC}

The following theorem removes the growth assumption for the
conclusion concerning infinitely many zeros. Let $n(r,a,f)$ count
the $a$-points of $f$ in $|z|\le r$, with multiplicity, and write
$n_p(r)=n(r,\infty,f)$ and $n_0(r)=n(r,0,f)$.

\begin{thm}[F.~Nazarov]\label{thm:finite-mass}
Let $(\zeta_j)_{j\ge1}$ be distinct complex numbers with
$|\zeta_j|\to\infty$, and let $(a_j)_{j\ge1}$ be nonzero complex
numbers satisfying \eqref{go}. Then
$f(z)=\sum_{j\ge1}a_j/(z-\zeta_j)$ has infinitely many zeros.
Moreover, for every $\eta>0$ there exists $r_\eta>0$ such that
\begin{equation}\label{eq:finite-mass-counts}
 n_p\left(\frac{r}{1+\eta}\right)-1
 \le n_0(r)\le n_p\bigl((1+\eta)r\bigr)-1
 \qquad(r\ge r_\eta).
\end{equation}
\end{thm}

A sketch of the proof of Theorem~\ref{thm:finite-mass} was explained
to the first-named author by Fedor Nazarov in 2015, but was not
published. We include the theorem with Nazarov's permission.
The comparison \eqref{eq:finite-mass-counts} does not by itself
establish the zero-deficiency conclusion of Theorem~\ref{thm:C}
without a growth assumption.

While this paper was in preparation, Shemyakov's
preprint~\cite[Theorem~1.1]{She26} appeared. It constructs a zero-free
Cauchy sum with positive residues, meromorphic order and lower order
$3/2$, and cumulative mass of order $r^{3/4}$. The construction
perturbs the zeros of a product of rotated Airy functions. Its
nonzero poles lie exactly on the rays $\arg z=\pm\pi/3$.
This gives a related example at potential order $3/4$, but does
not prescribe the residues to be the exact power sequence in
\eqref{eq:mass-profile}.

Our approach uses quasiconformal surgery and gives potentials
without critical points at every order $\rho\in[1/2,1)$, with the
mass sequences prescribed above. For $1/2<\rho<1$, the associated
Cauchy sums have order and lower order $\rho/(2\rho-1)$.
When $\rho=3/4$, our poles approach the rays
$\arg z=\pm\pi/3$ asymptotically.

\medskip
\noindent\textbf{Organization of the paper.}
The necessary tools are collected in Section~\ref{sec:tools}.
In Section~\ref{sec:boundary}, we extend a homeomorphism of the
imaginary axis to a locally quasiconformal homeomorphism of the
right half-plane. Section~\ref{sec:gluing} uses this extension to
deform the potential $\log|1+e^{-z}|$.
Section~\ref{sec:infinity} constructs the straightening map and
determines its behavior at infinity. In Section~\ref{sec:completion},
we complete the proofs of Theorems~\ref{thm:main} and~\ref{thm:endpoint}
and Corollary~\ref{cor:growth}. Finally, Section~\ref{sec:finite-mass}
proves Theorem~\ref{thm:finite-mass}.

\medskip
\noindent\textbf{Acknowledgments and AI tools disclosure.}
Teng Zhang is supported by the China Scholarship Council, the Young Elite
Scientists Sponsorship Program for PhD Students (China Association for Science
and Technology), and the Fundamental Research Funds for the Central
Universities at Xi'an Jiaotong University (Grant No.~xzy022024045).

ChatGPT 6.1 was used for language editing, proofreading, and grammatical
correction, and as an exploratory tool for discussing possible approaches to
selected results under the authors' mathematical supervision and guidance.
The authors take full responsibility for all mathematical arguments and for
the accuracy and correctness of the final manuscript.

\Needspace{7\baselineskip}
\section{Preliminaries}\label{sec:tools}

Let $\HH=\{w\in\C:\Rea w>0\}$ denote the right half-plane,
and $dA$ the planar Lebesgue measure. For an orientation-preserving
quasiconformal map $Q$, we write
$\mu_Q=Q_{\bar z}/Q_z$ for its Beltrami coefficient,
$DQ$ for its real differential, and $J_Q=\det DQ$ for
its Jacobian. The smaller singular value of $DQ$ is
denoted by $\smin(DQ)$. All statements involving
Beltrami coefficients are understood to hold almost
everywhere. For positive quantities $A$ and $B$,
the notation $A\asymp B$ means that their ratio is
bounded above and below by positive constants.

We first recall the measurable Riemann mapping theorem
which provides a normalized global solution of the Beltrami equation
and local uniqueness up to a conformal change of image coordinates.

\begin{lem}[{\cite[Chapter~V]{Ahl06}}]
\label{lem:mrmt}
Let $\mu$ be a measurable function on $\C$ satisfying
$\|\mu\|_\infty<1$. Then there exists a unique
quasiconformal homeomorphism $g:\C\to\C$ fixing
$0$, $1$ and $\infty$ such that
$g_{\bar z}=\mu g_z$ almost everywhere.
Moreover, if two locally quasiconformal homeomorphisms on the same
domain have the same Beltrami coefficient, then they differ by a
conformal change of image coordinates.
\end{lem}

To control the behavior of quasiconformal mappings at
infinity, we use the Teichm\"uller--Wittich--Belinskii
theorem in the following form. 

\begin{lem}[{\cite[Theorem~5.1 and Remark~5.2]{Bra10}}]
\label{lem:twb}
Let $G:\C\to\C$ be a quasiconformal homeomorphism
fixing $\infty$. Suppose that, for some $R>0$,
\begin{equation}\label{eq:twb}
 \int_{|w|>R}\frac{|\mu_G(w)|}{|w|^2}\dd A(w)<\infty.
\end{equation}
Then there exists a constant $c\in\C\setminus\{0\}$
such that $G(w)=cw+o(|w|)$ as $w\to\infty$,
uniformly in direction.
\end{lem}

We use the planar bi-Lipschitz Schoenflies theorem of
Tukia~\cite{Tuk80}.

\begin{lem}[{\cite[Theorem~A]{Tuk80}}]\label{lem:schoenflies}
Every bi-Lipschitz embedding of the unit circle into $\mathbb C$
extends to a bi-Lipschitz homeomorphism of $\mathbb C$.
\end{lem}

\begin{lem}\label{lem:harmonic-remainder}
Let $(\zeta_k)_{k\ge1}$ be a sequence of pairwise distinct
nonzero complex numbers with $|\zeta_k|\to\infty$, and let
$(a_k)_{k\ge1}$ be a sequence of positive numbers satisfying
\eqref{eq:sum}. Suppose that $v$ is subharmonic on $\C$,
harmonic outside the points $\zeta_k$, and that
$v(z)-a_k\log|z-\zeta_k|$ extends harmonically across
$\zeta_k$ for every $k\ge1$. If
\begin{equation}\label{eq:sublinear}
 \sup_{|z|=r}v^+(z)=o(r)\qquad(r\to\infty),
\end{equation}
where $v^+=\max\{v,0\}$, then
$v(z)-v(0)=\sum_{k=1}^{\infty}a_k
\log|1-z/\zeta_k|$.
\end{lem}

\begin{proof}
Put $U(z)=\sum_{k=1}^{\infty}a_k\log|1-z/\zeta_k|$.
By \eqref{eq:sum}, the series converges locally uniformly
away from the poles. For $|z|=r$, we have
$U^+(z)\le\sum_{k=1}^{\infty}
a_k\log(1+r/|\zeta_k|)$.
For each fixed $k$, the ratio
$r^{-1}\log(1+r/|\zeta_k|)$ tends to zero as $r\to\infty$
and is bounded by $1/|\zeta_k|$.
Since $\sum_{k=1}^{\infty}a_k/|\zeta_k|<\infty$,
dominated convergence gives
$\sup_{|z|=r}U^+(z)=o(r)$.

The function $h=v-U$ extends harmonically across
every $\zeta_k$ and is therefore harmonic on $\C$.
For a locally integrable function $\varphi$, let
$\langle\varphi\rangle_r$ denote its circular mean
on $|z|=r$. Since $v$ and $U$ are subharmonic,
the sub-mean inequality gives
$\langle v\rangle_r\ge v(0)$ and
$\langle U\rangle_r\ge U(0)=0$.
Using $|t|=2t^+-t$, we obtain
\[
 \langle |h|\rangle_r
 \le 2\langle v^++U^+\rangle_r-v(0)=o(r),
\]
where the last equality follows from
\eqref{eq:sublinear} and the estimate for $U^+$.

For each integer $m\ge1$, the coefficients of
$\cos(mt)$ and $\sin(mt)$ in the harmonic expansion
of $h(re^{it})$, divided by $r^m$, have absolute
value at most $2r^{-m}\langle|h|\rangle_r$.
Letting $r\to\infty$ shows that all nonconstant
coefficients vanish. Hence $h$ is constant.
Since $U(0)=0$, we have $h\equiv v(0)$, as required.
\end{proof}

For a meromorphic function $f$, let $m(r,f)$ denote the circular
mean of $\log^+|f|$ on $|z|=r$. When $f$ has no pole at the origin,
its Nevanlinna characteristic is
\[
 T(r,f)=m(r,f)+\int_0^r\frac{n_p(t)}t\dd t.
\]
We use the following proximity estimate for Cauchy sums, in the
form stated by Langley and Rossi~\cite[Theorem~1.1]{LR04};
see also~\cite[Chapter~V, Section~6]{GO08}.

\begin{lem}[{\cite[Theorem~1.1]{LR04}}]
\label{lem:proximity}
Let $(\zeta_k)_{k\ge1}$ be a sequence of pairwise
distinct complex numbers with $|\zeta_k|\to\infty$,
and let $(a_k)_{k\ge1}$ be a sequence of complex
numbers satisfying
$\sum_{\zeta_k\ne0}|a_k/\zeta_k|<\infty$.
Then the meromorphic function
$f(z)=\sum_{k=1}^{\infty}a_k/(z-\zeta_k)$
satisfies $m(r,f)=o(1)$ as $r\to\infty$.
\end{lem}

\section{Boundary data and their extension}\label{sec:boundary}

We construct a boundary homeomorphism for each of the three
sequences $(b_m)$ introduced above: the power sequence
in \eqref{eq:mass-profile} and the two endpoint sequences
defined by \eqref{eq:endpoint-models} and
\eqref{eq:endpoint-masses}. We refer to these as the power,
normal and maximal cases, respectively.

Recall that $L=2\pi$. In each case, define an odd, continuous,
piecewise affine function $\sigma:\R\to\R$ by linear
interpolation of the values
\begin{equation}\label{eq:sigma}
 \sigma(0)=0,\qquad
 \sigma(Lm)=L\sum_{\ell=1}^m b_\ell\quad(m\ge1).
\end{equation}
Since $b_m>0$ for every $m\ge1$ and
$\sum_{m\ge1}b_m=\infty$ in all three cases,
$\sigma$ is an increasing homeomorphism of $\R$.
For each $j\in\Z$, let $\alpha_j$ denote the slope
of $\sigma$ on $(Lj,L(j+1))$. By the definition
of $\sigma$ and its oddness,
\begin{equation}\label{eq:slopes}
 \alpha_{m-1}=\alpha_{-m}=b_m\quad(m\ge1).
\end{equation}
In all three cases, the slopes $\alpha_j$ are uniformly
bounded above and are bounded below by a positive constant
on every compact interval.

The following lemma gives the asymptotic estimates for
$\sigma$ needed to construct the quasiconformal maps.

\begin{lem}\label{lem:boundary}
In the power case, let $A_\beta=L^{1-\beta}/\beta$.
Then
\begin{equation}\label{eq:power-boundary}
 \sigma(r)\sim A_\beta r^\beta,\qquad
 0<c\le\frac{r\sigma'(r)}{\sigma(r)}\le C,\qquad
 \int_1^\infty
 \left|\frac{r\sigma'(r)}{\sigma(r)}-\beta\right|
 \frac{\dd r}{r}<\infty.
\end{equation}
In the normal and maximal cases, let $h$ be the
corresponding function in \eqref{eq:endpoint-models}.
Then
\begin{equation}\label{eq:endpoint-boundary}
 \sigma(r)=h(r)\left(1+O\left(\frac{|h''(r)|}{h(r)}\right)\right),
 \qquad
 \sigma'(r)=h'(r)(1+O(r^{-1})),\qquad
 \frac{|h''(r)|}{h'(r)}=O(r^{-1}).
\end{equation}
The estimates involving derivatives hold almost everywhere
for sufficiently large $r$.
\end{lem}

\begin{proof}
We first consider the power case. Put
$E(r)=\sigma(r)-A_\beta r^\beta$ for $r>0$.
For $L(m-1)<r<Lm$, where $m\ge2$, the slope
of $\sigma$ is $b_m$. Using \eqref{eq:mass-profile}
and the definition of $A_\beta$, we obtain
\begin{equation}\label{eq:boundary-error}
 |E'(r)|\le C r^{\beta-1}/m,\qquad
 \int_1^\infty |E'(r)|r^{-\beta}\dd r<\infty.
\end{equation}
For $r\ge1$, we have
$E(r)=E(1)+\int_1^r E'(t)\dd t$.
The integrability in \eqref{eq:boundary-error},
together with dominated convergence, gives
$E(r)/r^\beta\to0$ as $r\to\infty$.
Hence $\sigma(r)\sim A_\beta r^\beta$.

Applying Tonelli's theorem to the integral
representation of $E$, we also obtain
\begin{equation}\label{eq:boundary-integral}
 \int_1^\infty |E(r)|r^{-\beta-1}\dd r
 \le \frac{|E(1)|}{\beta}
       +\frac1\beta\int_1^\infty
       |E'(r)|r^{-\beta}\dd r<\infty.
\end{equation}
We have
$\sigma(r)\asymp r^\beta$ and
$\sigma'(r)\asymp r^{\beta-1}$ for sufficiently
large $r$. These estimates give the two-sided
bound for $r\sigma'(r)/\sigma(r)$.
Finally, the identity
$r\sigma'/\sigma-\beta=(rE'-\beta E)/\sigma$,
together with \eqref{eq:boundary-error} and
\eqref{eq:boundary-integral}, proves the integral
estimate in \eqref{eq:power-boundary}.

We now consider the normal and maximal cases.
By \eqref{eq:endpoint-masses}, we have
$\sigma(Lm)=h(Lm)$ for every $m\ge1$.
Thus $\sigma$ is the piecewise affine interpolation
of $h$ on the intervals $[L(m-1),Lm]$, $m\ge2$.
The standard interpolation estimate bounds
$|\sigma(r)-h(r)|$ by a constant times the supremum
of $|h''|$ on the interval containing $r$.

For either choice of $h$ in \eqref{eq:endpoint-models}, direct
differentiation gives $|h''(r)|/h'(r)=O(r^{-1})$.
Moreover, $h'$ and $|h''|$ vary by bounded factors
on intervals of the fixed length $L$.
The interpolation estimate therefore gives
$\sigma(r)=h(r)(1+O(|h''(r)|/h(r)))$.

Finally, the mean value theorem expresses the
slope of $\sigma$ on each interpolation interval
as $h'(s)$ for some point $s$ in that interval.
Since $|h''(r)|/h'(r)=O(r^{-1})$, we obtain
$\sigma'(r)=h'(r)(1+O(r^{-1}))$ almost everywhere.
This proves \eqref{eq:endpoint-boundary}.
\end{proof}

For $s>0$ and $t\in\R$, define
\begin{equation}\label{eq:ba}
 \begin{aligned}
 U(t,s)&=\frac1{2s}\int_{t-s}^{t+s}\sigma(v)\dd v,\\
 V(t,s)&=\frac1s\left(\int_t^{t+s}\sigma(v)\dd v
                         -\int_{t-s}^{t}\sigma(v)\dd v\right),\\
 T(s+it)&=V(t,s)+iU(t,s).
 \end{aligned}
\end{equation}
The map $T$ is the Beurling--Ahlfors extension
of the boundary homeomorphism $it\mapsto i\sigma(t)$,
with parameter $r=2$ in the notation of
Lehtinen~\cite[Section~1, formula~(2)]{Leh84}.
With this choice of parameter, affine boundary maps have conformal
extensions. In the normal and maximal cases, however, $\sigma$
need not be globally quasisymmetric. We therefore establish
the required properties of $T$ directly.

\begin{lem}\label{lem:extension}
The map $T$ defined by \eqref{eq:ba} is a locally
quasiconformal homeomorphism of $\HH$ onto itself.
It satisfies $T(\overline w)=\overline{T(w)}$ and extends
continuously to the boundary with $T(it)=i\sigma(t)$.
If $\sigma(v)=av+b$ for every $v\in[t-s,t+s]$, then
\begin{equation}\label{eq:affine}
 T(s+it)=a(s+it)+ib.
\end{equation}
Moreover,
\begin{equation}\label{eq:proper-extension}
 |T(w)|\asymp \sigma(|w|)
\end{equation}
uniformly as $|w|\to\infty$ in $\overline\HH$.

For every fixed $B>0$, the quasiconformal distortion of
$T$ is uniformly bounded on $0<\Rea w\le B$, and
$\Rea T(s+it)\le Cs$. Furthermore, there exists a bounded
nonnegative measurable function $\eta_B$ on $\R$ such that
\begin{equation}\label{eq:strip-envelope}
 |\mu_T(s+it)|\le\eta_B(y)
 \quad\text{if }0<s\le B,\ |t+y|\le\pi/2,\qquad
 \int_\R\frac{\eta_B(y)}{1+|y|}\dd y<\infty.
\end{equation}
\end{lem}

\begin{proof}
\emph{Step 1. Local invertibility.}
The boundary values, symmetry and \eqref{eq:affine}
follow directly from \eqref{eq:ba}.
Since $\sigma$ is strictly increasing, $V(t,s)>0$.
The global Lipschitz continuity of $\sigma$ gives
$V(t,s)\le Cs$, so $T$ maps $\HH$ into itself.

We first show that $T$ has positive Jacobian.
Fix $s>0$ and $t\in\R$, and put
$P=\sigma(t+s)-\sigma(t)$ and
$Q=\sigma(t)-\sigma(t-s)$.
Define
$u=s^{-1}\int_0^s(\sigma(t+v)-\sigma(t))\dd v$
and
$v_0=s^{-1}\int_0^s(\sigma(t)-\sigma(t-v))\dd v$.
Since $\sigma$ is strictly increasing,
$0<u<P$ and $0<v_0<Q$.
Differentiating the expressions in \eqref{eq:ba} gives
\begin{equation}\label{eq:ba-jacobian}
 J_T(s+it)=\frac{2PQ-Qu-Pv_0}{s^2}>0.
\end{equation}
Indeed, the numerator equals
$Q(P-u)+P(Q-v_0)$, which is strictly positive.

The map $T$ is continuously differentiable in $\HH$.
Suppose that the slopes of $\sigma$ on
$[t-s,t+s]$ lie between $c\ell$ and $C\ell$,
where $\ell>0$. Then $P,Q\ge c\ell s$ and
$P-u,Q-v_0\ge c\ell s/2$.
Together with the derivatives of the averages
in \eqref{eq:ba}, these bounds yield
\begin{equation}\label{eq:local-differential}
 \|DT\|\le C_1\ell,\qquad J_T\ge c_1\ell^2,\qquad
 \smin(DT)\ge c_2\ell.
\end{equation}
Thus $T$ is a local diffeomorphism with locally bounded
quasiconformal distortion. Since the slopes of $\sigma$
are bounded above and below by positive constants on
compact intervals, the distortion is also uniformly
bounded in a neighborhood of every compact boundary
interval.

\emph{Step 2. Properness.}
We next prove \eqref{eq:proper-extension}.
By Lemma~\ref{lem:boundary}, the restriction of
$\sigma$ to the positive real axis is regularly varying
of index $\beta$ in the power case and of index zero
in the normal and maximal cases.

First suppose that $|t|\ge2s$ and $|t|\to\infty$.
The interval $[t-s,t+s]$ lies entirely on one side
of the origin, and $|v|\asymp|t|$ throughout
this interval. Regular variation gives
$|U(t,s)|\asymp\sigma(|t|)$ and
$V(t,s)=O(\sigma(|t|))$.
Since $|s+it|\asymp|t|$, we obtain
$|T(s+it)|\asymp\sigma(|s+it|)$.

It remains to consider $|t|\le2s$, where
$s\to\infty$ as $|s+it|\to\infty$.
Put $\tau=t/s$, so that $|\tau|\le2$.
The regular variation of $\sigma$, together with
its oddness, gives
$\sigma(sv)/\sigma(s)\to\sgn(v)|v|^\beta$
in the power case and
$\sigma(sv)/\sigma(s)\to\sgn(v)$
in the normal and maximal cases, for every $v\ne0$.
The bound $\sigma(3s)/\sigma(s)=O(1)$ allows
dominated convergence on $[-3,3]$.
Hence the convergence also holds in $L^1[-3,3]$.

It follows that the normalized integrals defining
$U$ and $V$ converge uniformly for $|\tau|\le2$.
In the power case, the limit value of
$V(t,s)/\sigma(s)$ is strictly positive for every
$\tau\in[-2,2]$, since $v\mapsto\sgn(v)|v|^\beta$
is strictly increasing.

In the normal and maximal cases, direct integration
of the limit function $\sgn(v)$ gives
\begin{equation}\label{eq:slow-limit}
 \frac{U(t,s)}{\sigma(s)}
    =\max\{-1,\min\{\tau,1\}\}+o(1),\qquad
 \frac{V(t,s)}{\sigma(s)}
    =2(1-|\tau|)^++o(1).
\end{equation}
The errors are uniform for $|\tau|\le2$.
The two limit components in \eqref{eq:slow-limit}
never vanish simultaneously. Consequently,
$|T(s+it)|\asymp\sigma(s)$ in both endpoint cases.
Since $|s+it|\asymp s$ when $|t|\le2s$,
this proves \eqref{eq:proper-extension} in all
three cases.

We now show that $T:\HH\to\HH$ is proper.
By \eqref{eq:proper-extension}, the inverse image
of a compact subset of $\HH$ cannot contain a
sequence tending to infinity. Nor can such a
sequence approach a finite boundary point, because
$0<\Rea T(s+it)\le Cs$.
Thus the inverse image of every compact subset
of $\HH$ is compact in $\HH$.

By \eqref{eq:ba-jacobian}, $T$ is a local
diffeomorphism. Its image is open, and properness
implies that the image is also closed in $\HH$.
Since $\HH$ is connected, $T$ is surjective.
Indeed, properness and local injectivity imply that each
fiber is finite. Around the finitely many points of a fiber,
choose pairwise disjoint neighborhoods on which $T$ is a
diffeomorphism. Properness allows their image neighborhoods
to be shrunk to a common neighborhood with no additional
preimages. Thus $T$ is a covering map. Since its domain is
connected and $\HH$ is simply connected, the covering
has one sheet, so $T$ is a homeomorphism of $\HH$
onto itself.

\emph{Step 3. Estimates in a strip.}
Fix $B>0$. We prove the uniform distortion
bound and \eqref{eq:strip-envelope} on
$0<s\le B$.

Consider first the power case. Suppose that
$|t-Lm|\le L$, where $m$ is sufficiently large.
Since $s\le B$, only finitely many slopes of
$\sigma$ occur on $[t-s,t+s]$.
More precisely, their indices $j$ satisfy
$|j-m|\le N_B$ for some integer $N_B$
depending only on $B$.
Put $\ell_m=m^{\beta-1}$.
By \eqref{eq:mass-profile}, these slopes are
uniformly comparable to $\ell_m$, and
\begin{equation}\label{eq:strip-power}
 \frac{|\sigma'(v)-\ell_m|}{\ell_m}
 \le C_B m^{-1}.
\end{equation}
This estimate holds almost everywhere on the
intervals under consideration.

Subtract from $\sigma$ the affine function
with slope $\ell_m$ that agrees with $\sigma$
at $t$. The extension of this affine function
is conformal by \eqref{eq:affine}.
Since the operators defining $U$ and $V$
in \eqref{eq:ba} are linear, \eqref{eq:strip-power}
bounds the differential of the difference
between the two extensions by $\ell_m$ times
the right-hand side of \eqref{eq:strip-power}.
Combining this estimate with the lower bound
for $\smin(DT)$ in \eqref{eq:local-differential},
we obtain the corresponding bound for $|\mu_T|$.
The same differential estimates also give a
uniform bound on the quasiconformal distortion.
For large negative $t$, the estimates follow
from the oddness of $\sigma$ and the symmetry of $T$.

In the normal and maximal cases,
\eqref{eq:endpoint-boundary} gives
$\sigma'(v)=h'(|t|)(1+O_B(|t|^{-1}))$
whenever $|v-t|\le B$ and $|t|$ is sufficiently
large. Subtracting the affine function of slope
$h'(|t|)$ and applying \eqref{eq:affine} as above,
we obtain
\begin{equation}\label{eq:strip-endpoint}
 DT(s+it)=h'(|t|)(I+O_B(|t|^{-1})),\qquad
 |\mu_T(s+it)|\le C_B/(1+|t|).
\end{equation}
Here $I$ denotes the $2\times2$ identity matrix.
Together with the local bounds from Step~1,
these estimates prove that the distortion of
$T$ is uniformly bounded on every strip
$0<s\le B$.

It remains to construct $\eta_B$.
In the power case, enlarge the index range
in \eqref{eq:strip-power} by a fixed number
of indices to account for the additional
restriction $|t+y|\le\pi/2$.
On each corresponding interval of the
$y$-axis, choose $\eta_B(y)$ to be a sufficiently
large constant multiple of the enlarged
right-hand side of \eqref{eq:strip-power},
truncated at $1$. Set $\eta_B(y)=1$ on
a bounded interval containing the origin.
Then we obtain
\[
 \int_\R\frac{\eta_B(y)}{1+|y|}\dd y
 \le C_B\left(1+
 \sum_{m\ge1}m^{-2}\right)
 <\infty,
\]
where the convergence follows from
\eqref{eq:mass-profile}.

In the normal and maximal cases,
\eqref{eq:strip-endpoint} allows us to take
$\eta_B(y)=\min\{1,C_B/(1+|y|)\}$,
with $C_B$ sufficiently large to cover the
bounded values of $y$ as well.
This function satisfies the required integral
condition and completes the proof of
\eqref{eq:strip-envelope}.
\end{proof}

\section{A quasiconformal deformation of the potential}\label{sec:gluing}

We begin with the logarithmic potential
\begin{equation}\label{eq:base}
 u_0(z)=\log|1+e^{-z}|,\qquad
 (u_0)_x-i(u_0)_y=-\frac1{1+e^z}.
\end{equation}
The conjugate gradient of $u_0$ has no zeros.
The logarithmic poles of $u_0$ are
$\xi_j=(2j+1)\pi i$, $j\in\Z$.
The connected components of the set $\{u_0<0\}$ are
\begin{equation}\label{eq:negative-components}
 V_j=\{x+iy:\pi/2+Lj<y<3\pi/2+Lj,\
                         x>-\log(-2\cos y)\}.
\end{equation}
Each component $V_j$ contains exactly one logarithmic
pole, namely $\xi_j$.

Let $D=\{z\in\C:u_0(z)>0\}$.
The domain $D$ is connected. Indeed, every point of $D$
can be joined to the half-plane $\{\Rea z<-\log2\}$
by a horizontal segment contained in $D$.

\begin{lem}\label{lem:coordinate}
There exists a conformal map $H:D\to\HH$ that takes
real values on the real axis and satisfies
\begin{equation}\label{eq:coordinate}
 e^{H(z)}=1+e^{-z},\qquad \Rea H(z)=u_0(z),\qquad
 H(z+Li)=H(z)-Li.
\end{equation}
The map $H$ extends conformally across every finite
boundary point of $D$. For each $j\in\Z$, its boundary
extension maps $\partial V_j$ onto
$i(-L(j+1),-Lj)$.

Moreover, for every fixed $B>0$,
\begin{equation}\label{eq:coordinate-bounds}
 |\Ima H(z)+\Ima z|<\pi/2,\qquad
 \Rea H(z)\le\log(1+e^B)\quad(\Rea z\ge-B,\ z\in D).
\end{equation}
On the half-plane $\Rea z<-\log2$, we have
\begin{equation}\label{eq:left-coordinate}
 H(z)=-z+\Log(1+e^z)\qquad(\Rea z<-\log2).
\end{equation}
\end{lem}

\begin{proof}
For $w\in\HH$, define
\begin{equation}\label{eq:inverse-coordinate}
 h_0(w)=-w-\Log(1-e^{-w}).
\end{equation}
Since $|e^{-w}|<1$, the logarithm is analytic
on $\HH$ and is given by its convergent power
series. Moreover, $1-e^{-w}$ has positive
real part.

Direct calculation gives
$1+e^{-h_0(w)}=e^w$ and
$h_0(w+Li)=h_0(w)-Li$.
We claim that $h_0$ is injective.
If $h_0(w_1)=h_0(w_2)$, then
$e^{w_1}=e^{w_2}$, so
$w_1-w_2=Lmi$ for some $m\in\Z$.
The translation identity then gives
$h_0(w_1)=h_0(w_2)-Lmi$, and hence $m=0$.

To prove surjectivity onto $D$, fix $z\in D$
and choose a logarithm $w$ of $1+e^{-z}$.
Since $u_0(z)>0$, we have $\Rea w>0$.
The identity $1+e^{-h_0(w)}=e^w$
implies $h_0(w)=z+Lmi$ for some $m\in\Z$.
Consequently, $h_0(w+Lmi)=z$.
Thus $h_0$ maps $\HH$ conformally onto $D$,
and $H=h_0^{-1}$.

The identities in \eqref{eq:coordinate} follow
from the definition of $h_0$. Since $h_0$
maps the positive real axis onto the real axis,
$H$ is real on $\R$.

We next examine the boundary correspondence.
For $0<v<L$, the boundary values of $h_0$
are
\begin{equation}\label{eq:boundary-coordinate}
 h_0(iv)=-\log(2\sin(v/2))-i(v/2+\pi/2).
\end{equation}
As $v$ runs from $0$ to $L$, the expression
in \eqref{eq:boundary-coordinate} traces
$\partial V_{-1}$ exactly once, with nonzero
derivative. The translation identity
$h_0(w+Li)=h_0(w)-Li$ gives the boundary
correspondence for every $\partial V_j$.

At a finite point of $\partial V_j$,
the function $1+e^{-z}$ is nonzero and has
a local analytic logarithm. Its derivative is
$-1/(1+e^z)$, which does not vanish at such
a point. Hence $H$ extends conformally
across every finite boundary point.
No finite boundary point is mapped to a node
$iLj$, since $e^{H(z)}=1$ would imply
$e^{-z}=0$.

Figure~\ref{fig:coordinate} illustrates the boundary
correspondence for $V_0$ and $V_{-1}$.
In particular, $H$ maps $\partial V_0$ onto
$i(-L,0)$ and $\partial V_{-1}$ onto $i(0,L)$.

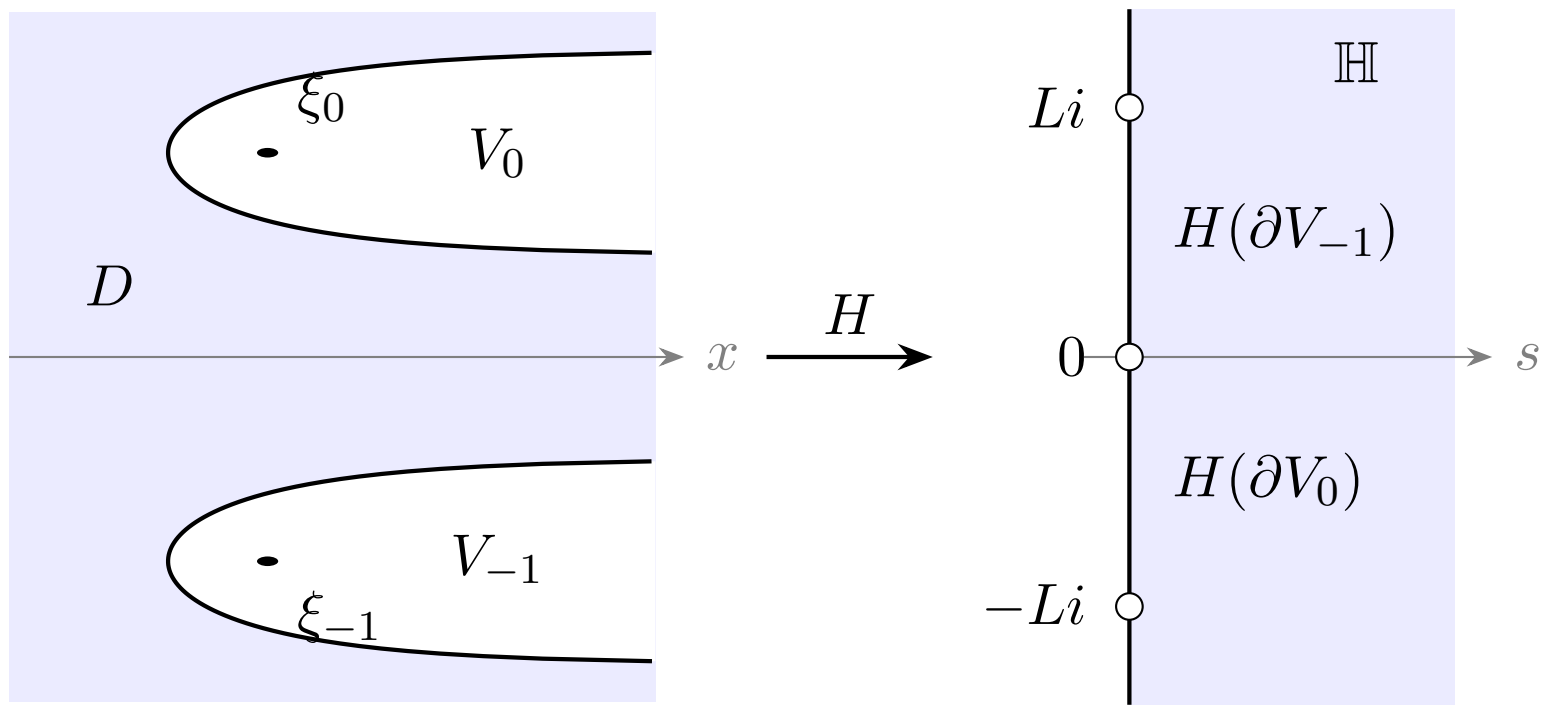
\begin{figure}[htbp]
\centering
\begin{tikzpicture}[>=Stealth,font=\small]
 \begin{scope}[xscale=.95,yscale=.43]
  \fill[blue!8] (-1.8,-5.3) rectangle (2.7,5.3);
  \foreach \c in {-3.14159,3.14159}{
   \fill[white] plot[domain=-1.537:1.537,samples=80]
     ({-ln(2*cos(\x r))},{\c+\x}) -- (2.7,{\c+1.537}) -- (2.7,{\c-1.537}) -- cycle;
   \draw[thick] plot[domain=-1.537:1.537,samples=80]
     ({-ln(2*cos(\x r))},{\c+\x});
   \fill (0,\c) circle[radius=2.1pt];
  }
  \draw[->,gray] (-1.8,0)--(2.9,0) node[right]{$x$};
  \node at (-1.1,1.1) {$D$};
  \node at (1.6,3.14159) {$V_0$};
  \node at (1.6,-3.14159) {$V_{-1}$};
  \node[above right=2pt] at (0,3.14159) {$\xi_0$};
  \node[below right=2pt] at (0,-3.14159) {$\xi_{-1}$};
 \end{scope}
 \draw[->,thick] (3.3,0)--(4.4,0) node[midway,above]{$H$};
 \begin{scope}[shift={(5.7,0)}]
  \fill[blue!8] (0,-2.3) rectangle (2.15,2.3);
  \draw[thick] (0,-2.3)--(0,2.3);
  \draw[->,gray] (-.3,0)--(2.4,0) node[right]{$s$};
  \node at (1.5,1.95) {$\HH$};
  \foreach \y/\lab in {-1.65/{-Li},0/{0},1.65/{Li}}{
   \draw[fill=white] (0,\y) circle[radius=2.5pt];
   \node[left=4pt] at (0,\y) {$\lab$};
  }
  \node[right=4pt] at (0,-.825) {$H(\partial V_0)$};
  \node[right=4pt] at (0,.825) {$H(\partial V_{-1})$};
 \end{scope}
\end{tikzpicture}
\caption{The conformal map $H$ sends $D$ onto the right
half-plane. The boundary curves $\partial V_0$ and
$\partial V_{-1}$ correspond to the open intervals
$i(-L,0)$ and $i(0,L)$, respectively. The endpoints
of these intervals correspond to the infinite ends
of the boundary curves.}
\label{fig:coordinate}
\end{figure}

It remains to prove the estimates.
For $w\in\HH$, the number $1-e^{-w}$
lies in the right half-plane. Taking imaginary
parts in \eqref{eq:inverse-coordinate} therefore
gives $|\Ima h_0(w)+\Ima w|<\pi/2$.
Substituting $w=H(z)$ proves the first
inequality in \eqref{eq:coordinate-bounds}.
For $\Rea z\ge-B$, we have
$|1+e^{-z}|\le1+e^B$, and hence
$\Rea H(z)\le\log(1+e^B)$.

Finally, when $\Rea z<-\log2$, the function
$\Log(1+e^z)$ is analytic and is defined
by its power series. The expression
$-z+\Log(1+e^z)$ satisfies
$e^{H(z)}=1+e^{-z}$ and agrees with $H$
on the real axis. The identity theorem
then gives \eqref{eq:left-coordinate}.
\end{proof}

We now combine the conformal map $H$ with the
Beurling--Ahlfors extension $T$ constructed in
Lemma~\ref{lem:extension}. Define a measurable
Beltrami coefficient $\mu$ and an extended
real-valued function $\widetilde u$ by
\begin{equation}\label{eq:physical-mu}
 \mu(z)=
 \begin{cases}
  \mu_T(H(z))\,\overline{H'(z)}/H'(z),&z\in D,\\
  0,&z\notin D,
 \end{cases}
 \qquad
 \widetilde u(z)=
 \begin{cases}
  \Rea T(H(z)),&z\in D,\\
  \alpha_j u_0(z),&z\in V_j,\\
  0,&u_0(z)=0.
 \end{cases}
\end{equation}
At each logarithmic pole $\xi_j$, we set
$\widetilde u(\xi_j)=-\infty$.

We verify that the definitions in
\eqref{eq:physical-mu} agree across the boundary
curves. Fix $j\in\Z$ and a finite point
$z_0\in\partial V_j$. By Lemma~\ref{lem:coordinate},
the boundary value $H(z_0)$ belongs to the interior
of $i(-L(j+1),-Lj)$.
On this interval, $\sigma$ is affine with slope
$\alpha_{-j-1}=\alpha_j$, by \eqref{eq:slopes}.
Thus \eqref{eq:affine} gives
$\Rea T(H(z))=\alpha_j\Rea H(z)=\alpha_j u_0(z)$
for $z\in D$ sufficiently close to $z_0$.
The same affine identity shows that
$\mu_T(H(z))=0$ there.
Consequently, $\widetilde u=\alpha_j u_0$
and $\mu=0$ throughout a neighborhood of $z_0$.
In particular, the functions defining
$\widetilde u$ agree harmonically across
$\partial V_j$.

By Lemma~\ref{lem:extension}, the coefficient
$\mu$ is locally bounded away from modulus one.
For each fixed $B>0$, this bound is uniform on
the half-plane $\Rea z\ge-B$.
We also claim that
\begin{equation}\label{eq:physical-integrability}
 \int_{\substack{x\ge0\\|z|>2}}
          \frac{|\mu(z)|}{|z|^2}\dd A(z)<\infty.
\end{equation}
To prove this, write $z=x+iy$ with $x\ge0$.
For $z\in D$, \eqref{eq:coordinate-bounds} gives
$0<\Rea H(z)\le\log2$ and
$|\Ima H(z)+y|<\pi/2$.
Applying \eqref{eq:strip-envelope} with
$B=\log2$, we obtain
$|\mu(z)|\le\eta_{\log2}(y)$.
Since $\mu=0$ outside $D$, this estimate
holds throughout the right half-plane.

For $|y|\ge1$, integration with respect to
$x$ gives
$\int_0^\infty(x^2+y^2)^{-1}\dd x
=\pi/(2|y|)$.
The integral over $|y|\ge1$ is therefore
finite by the weighted integrability of
$\eta_{\log2}$ in \eqref{eq:strip-envelope}.
For $|y|<1$, the restriction $|z|>2$
keeps $x$ bounded away from zero.
Since $|\mu|\le1$, the remaining integral
is finite as well. This proves
\eqref{eq:physical-integrability}.

The following lemma shows that a locally quasiconformal
homeomorphism with Beltrami coefficient $\mu$
transforms $\widetilde u$ into a subharmonic
function with no critical points away from its
logarithmic poles. The construction of such a
global homeomorphism is given in
Section~\ref{sec:infinity}.

\begin{lem}\label{lem:gluing}
Let $g:\C\to\C$ be an orientation-preserving,
locally quasiconformal homeomorphism with
Beltrami coefficient $\mu$ defined by
\eqref{eq:physical-mu}.
Then $v=\widetilde u\circ g^{-1}$ is subharmonic
on $\C$ and harmonic with nonvanishing gradient
on $\C\setminus\{g(\xi_j):j\in\Z\}$.
Moreover, for every $j\in\Z$, the function
$v(Z)-\alpha_j\log|Z-g(\xi_j)|$ extends
harmonically across $g(\xi_j)$.
\end{lem}

\begin{proof}
On $D$, the maps $T\circ H$ and $g$ have the
same Beltrami coefficient $\mu$.
By Lemma~\ref{lem:mrmt}, the map
$T\circ H\circ g^{-1}$ is conformal on $g(D)$.
Consequently,
$v=\Rea(T\circ H\circ g^{-1})$ is harmonic
on $g(D)$ and has nonvanishing gradient.

On each component $V_j$, the coefficient
$\mu$ vanishes, so $g$ is conformal.
Since $\widetilde u=\alpha_j u_0$ on $V_j$,
we have
$v=\alpha_j u_0\circ g^{-1}$ on $g(V_j)$.
The derivative formula in \eqref{eq:base}
shows that $u_0$ has nonvanishing gradient
away from $\xi_j$. Thus $v$ is harmonic
with nonvanishing gradient on
$g(V_j)\setminus\{g(\xi_j)\}$.

Now let $z_0$ be a finite point of
$\partial V_j$. As established after
\eqref{eq:physical-mu}, there is a
neighborhood of $z_0$ on which
$\mu=0$ and $\widetilde u=\alpha_j u_0$.
The map $g$ is conformal on this
neighborhood, so
$v=\alpha_j u_0\circ g^{-1}$
is harmonic with nonvanishing gradient
across the image of the boundary curve.

Figure~\ref{fig:gluing-mechanism} illustrates this matching in a
neighborhood of the seam and the preservation of the nonvanishing
gradient under the conformal change of coordinates there.

\begin{figure}[htbp]
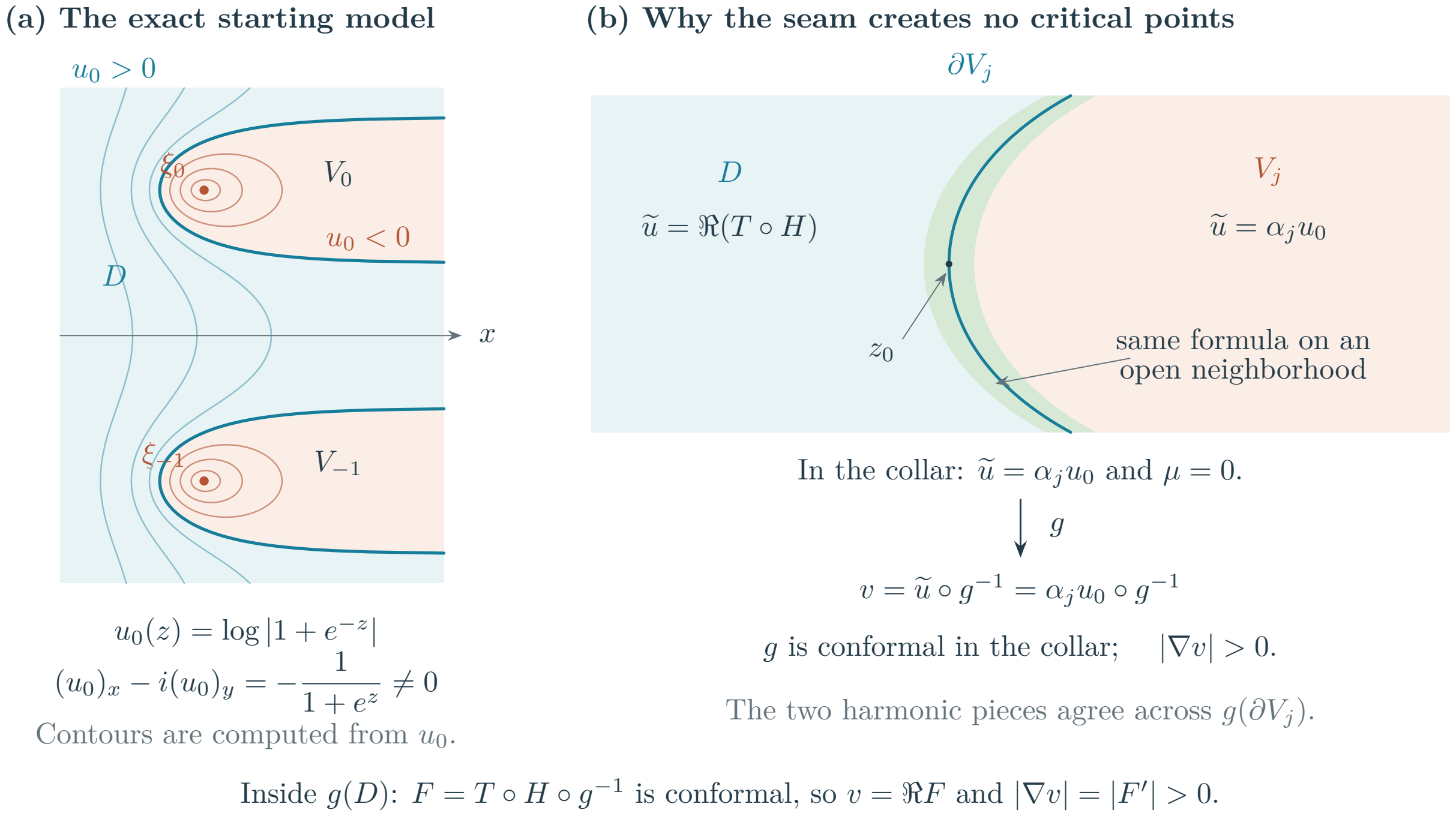

\centering
\begingroup
\resizebox{\linewidth}{!}{%
%
}
\endgroup
\caption{The local mechanism in Lemma~\ref{lem:gluing}.
(a) Level sets of the starting potential $u_0(z)=\log|1+e^{-z}|$;
each negative component contains one logarithmic pole.
(b) Near every finite point of $\partial V_j$, both definitions of
$\widetilde u$ agree with $\alpha_j u_0$ and the Beltrami coefficient
vanishes. Thus $g$ is conformal on this neighborhood, and the two
pieces of $v$ join with nonvanishing gradient. The width of the
highlighted neighborhood is schematic.}
\label{fig:gluing-mechanism}
\end{figure}

Finally, near $\xi_j$, the function
$u_0(z)-\log|z-\xi_j|$ extends harmonically
across $\xi_j$. Since $g$ is conformal
near $\xi_j$, it follows that
$v(Z)-\alpha_j\log|Z-g(\xi_j)|$
extends harmonically across $g(\xi_j)$.
The logarithmic coefficient at
$g(\xi_j)$ is therefore $\alpha_j>0$.
Hence $v$ is subharmonic on $\C$,
with logarithmic poles at $g(\xi_j)$
and no critical points elsewhere.
\end{proof}

\section{A comparison homeomorphism and asymptotics of \texorpdfstring{$g$}{g}}
\label{sec:infinity}

We construct a comparison homeomorphism $\Phi$ whose
Beltrami coefficient agrees with $\mu$ sufficiently far
to the left. The remaining distortion, expressed in the
$\Phi$-coordinate, will satisfy the integral condition
\eqref{eq:twb}. This approach is particularly important
in the normal and maximal cases, where $\mu$ need not
be uniformly bounded away from modulus one on $\C$.

\subsection{The comparison domains}

In the power case, put
$q=1/\rho=2/(1+\beta)$, so that $p=2-q=\beta q$.
For the power, normal and maximal cases, choose the
constant $A$ and the conformal map $\Psi$ according to
\begin{equation}\label{eq:comparison-functions}
 \begin{array}{c|c|c}
  \text{case}&A&\Psi(w)\\ \hline
  \text{power}&0&-w^q\\
  \text{normal}&\pi/2&-w^2\\
  \text{maximal}&A\text{ sufficiently large}&-w^2/\Log w .
 \end{array}
\end{equation}
All branches are chosen on the right half-plane,
with $\Log w$ real for $w>0$. Write
$\HH_A=\{w\in\C:\Rea w>A\}$.
For each choice in \eqref{eq:comparison-functions},
the function $\Psi$ is conformal on $\HH_A$.

In the power and normal cases, this follows directly
from the choice of branches. For the maximal case,
define $k(w)=w/\sqrt{\Log w}$ on $\HH_A$.
When $A$ is sufficiently large, $\Rea k'(w)>0$
throughout $\HH_A$. Since $\HH_A$ is convex,
integration of $k'$ along the line segment joining
two points of $\HH_A$ shows that $k$ is injective.
Moreover, $\arg k(w)$ has the same sign as $\arg w$
and smaller absolute value, so $k(\HH_A)\subset\HH$.
Squaring is injective on $\HH$. Hence
$\Psi(w)=-k(w)^2$ is conformal on $\HH_A$.

We first define a comparison map $\Phi_0$ outside
a sufficiently large disk. For $\Rea z\le0$, put
\begin{equation}\label{eq:phi-left}
 \Phi_0(z)=\Psi(A+T(-z)).
\end{equation}
For $\Rea z\ge0$, write $z=re^{it}$ with
$|t|\le\pi/2$ and define
\begin{equation}\label{eq:phi-right}
 \Phi_0(re^{it})=
 \begin{cases}
  \sigma(r)^q e^{ipt},&\text{power case},\\
  (\sigma(r)+it)^2,&\text{normal case},\\
  \mathcal R(r)\exp(2it\omega(r)/\pi),&\text{maximal case},
 \end{cases}
\end{equation}
where, in the maximal case,
\begin{equation}\label{eq:polar-data}
 \mathcal R(r)=|\Psi(A-i\sigma(r))|,\qquad
 \omega(r)=\arg\Psi(A-i\sigma(r))\in(0,\pi)
\end{equation}
for sufficiently large $r$.

The definitions in \eqref{eq:phi-left} and
\eqref{eq:phi-right} agree on the imaginary axis.
Indeed, for $r>0$, the boundary values of $T$
give $T(-ir)=-i\sigma(r)$.
In the power case, the identity $p=2-q$ gives
$-(-i\sigma(r))^q=\sigma(r)^q e^{ip\pi/2}$.
In the normal case, $A=\pi/2$ gives
$-(A-i\sigma(r))^2=(\sigma(r)+i\pi/2)^2$.
In the maximal case, the two boundary values
agree by the definitions of $\mathcal R$ and
$\omega$ in \eqref{eq:polar-data}.
The matching on the negative imaginary axis
follows by complex conjugation.

\begin{lem}\label{lem:comparison}
The map $\Phi_0$ defined by
\eqref{eq:phi-left}--\eqref{eq:polar-data}
extends across a bounded set to a locally
quasiconformal homeomorphism of $\C$.
The extension commutes with complex conjugation,
fixes zero and preserves the upper half-plane.

For every fixed $B>0$, the quasiconformal
distortion of $\Phi_0$ is uniformly bounded
on $\Rea z\ge-B$. Moreover, for sufficiently
large $|z|$ in this half-plane,
\begin{equation}\label{eq:comparison-jacobian}
 \frac{J_{\Phi_0}(z)}{|\Phi_0(z)|^2}
 \le \frac{C_B}{|z|^2}.
\end{equation}
On the right half-plane, $\Phi_0$ satisfies
the stronger integral estimate
\begin{equation}\label{eq:comparison-integral}
 \int_{\substack{x\ge0\\|z|>R}}
 |\mu_{\Phi_0}(z)|\,
 \frac{J_{\Phi_0}(z)}{|\Phi_0(z)|^2}\dd A(z)<\infty
\end{equation}
for all sufficiently large $R$.
\end{lem}

\begin{proof}
\emph{Step 1. Estimates on the right half-plane.}
Write $z=re^{it}$, where $r$ is sufficiently
large and $|t|\le\pi/2$. We use logarithmic
polar coordinates $s=\log r$ and $t=\arg z$.

In the power case, \eqref{eq:phi-right} gives
$\log|\Phi_0(re^{it})|=q\log\sigma(r)$
and $\arg\Phi_0(re^{it})=pt$.
Thus the differential in logarithmic polar
coordinates is
$\operatorname{diag}(qr\sigma'(r)/\sigma(r),p)$.
By \eqref{eq:power-boundary}, both diagonal
entries are bounded above and below by
positive constants. Since $q\beta=p$,
the same estimates give
\begin{equation}\label{eq:power-comparison}
 |\mu_{\Phi_0}(z)|
       \le C\left|\frac{r\sigma'(r)}{\sigma(r)}
                   -\beta\right|,
 \qquad
 \frac{J_{\Phi_0}(z)}{|\Phi_0(z)|^2}
       \le C/r^2 .
\end{equation}
Integrating the product of these two bounds
in logarithmic polar coordinates and using
\eqref{eq:power-boundary}, we obtain
\eqref{eq:comparison-integral} in the power case.

For the normal case, \eqref{eq:endpoint-boundary}
with $h(r)=\log r$ gives
$\sigma(e^s)=s+O(e^{-2s})$ and
$d\sigma(e^s)/ds=1+O(e^{-s})$.
By \eqref{eq:phi-right}, the map $\Phi_0$
is the composition of
$(s,t)\mapsto\sigma(e^s)+it$
with the conformal squaring map.
The first map has uniformly bounded distortion
for large $s$. Since $\sigma(r)\sim\log r$
and $|t|\le\pi/2$, direct calculation gives
\begin{equation}\label{eq:normal-comparison}
 \frac{J_{\Phi_0}(z)}{|\Phi_0(z)|^2}
       \le\frac{C}{r^2\log^2 r}.
\end{equation}
The right-hand side is integrable over
$\{\Rea z\ge0,\ |z|>R\}$.
As $|\mu_{\Phi_0}|<1$ almost everywhere,
\eqref{eq:normal-comparison} proves
\eqref{eq:comparison-integral} in the normal case.

We now consider the maximal case.
Put $\ell=\log\sigma(e^s)$.
The interpolation estimates in
\eqref{eq:endpoint-boundary}, applied to
$h(r)=\exp\sqrt{\log r}$, give
\begin{equation}\label{eq:maximal-interpolation}
 \ell=\sqrt s+O(e^{-2s}/\sqrt s),\qquad
 \frac{\dd\ell}{\dd s}
       =\frac1{2\sqrt s}(1+O(e^{-s})).
\end{equation}
Using $\Psi(w)=-w^2/\Log w$ and
\eqref{eq:polar-data}, we obtain
\begin{equation}\label{eq:polar-expansion}
 \begin{aligned}
  \log\mathcal R(e^s)&=2\ell-\log\ell+O(\ell^{-2}),\\
  \omega(e^s)&=\frac{\pi}{2\ell}
                   +O(\ell^{-3})+O(e^{-\ell}),\\
  \frac{\dd\log\mathcal R(e^s)}{\dd s}
       &=s^{-1/2}-\tfrac12s^{-1}+O(s^{-2}),\\
  \frac{2\omega(e^s)}{\pi}
       &=s^{-1/2}+O(s^{-3/2}),\\
  \frac{\dd\omega(e^s)}{\dd s}
       &=-\frac{\pi}{4}s^{-3/2}+O(s^{-5/2}).
 \end{aligned}
\end{equation}
The derivative estimates follow by differentiating
the expressions in \eqref{eq:polar-data} and
applying \eqref{eq:maximal-interpolation}.

In logarithmic polar coordinates,
\eqref{eq:phi-right} sends $(s,t)$ to
$(\log\mathcal R(e^s),2t\omega(e^s)/\pi)$.
Its differential is therefore
\[
 \begin{pmatrix}
  \dd\log\mathcal R(e^s)/\dd s&0\\
  (2t/\pi)\dd\omega(e^s)/\dd s&2\omega(e^s)/\pi
 \end{pmatrix}
 =s^{-1/2}(I+O(s^{-1/2})).
\]
The error estimate is uniform for
$|t|\le\pi/2$.
In particular, $\Phi_0$ has uniformly
bounded distortion for large $r$.
The same matrix estimate yields
\begin{equation}\label{eq:maximal-comparison}
 |\mu_{\Phi_0}(z)|
       \le C/\sqrt{\log r},\qquad
 \frac{J_{\Phi_0}(z)}{|\Phi_0(z)|^2}
       \le\frac{C}{r^2\log r}.
\end{equation}
Consequently,
$|\mu_{\Phi_0}(z)|J_{\Phi_0}(z)/
|\Phi_0(z)|^2$
is bounded by a constant multiple of
$r^{-2}(\log r)^{-3/2}$.
This function is integrable over the
right half-plane outside a sufficiently
large disk, proving
\eqref{eq:comparison-integral} in the
maximal case.

Thus the integral estimate holds in
all three cases. The estimates above
also establish the bounded distortion
of $\Phi_0$ on the right half-plane.

\emph{Step 2. Injectivity and extension across a compact set.}
We first show that the formulas
\eqref{eq:phi-left} and \eqref{eq:phi-right}
define a homeomorphism outside a sufficiently
large disk.

In the power case, $\Psi$ maps $\HH$
onto a sector of angle $q\pi$.
The image of the left half-plane under
\eqref{eq:phi-left} lies in this sector,
while \eqref{eq:phi-right} maps the right
half-plane into the complementary sector
of angle $p\pi$. Both maps are injective
on their respective domains, and their
boundary values agree on the imaginary axis.

In the normal case, the common boundary
of the two image domains is the parabola
$\Psi(A+i\R)$, with $A=\pi/2$.
The left formula maps into $\Psi(\HH_A)$,
while the right formula maps into the
opposite side of the parabola.
To see the injectivity of the right formula,
take the square root with positive real part.
Equation \eqref{eq:phi-right} then gives
$\sqrt{\Phi_0(re^{it})}=\sigma(r)+it$.
Since $\sigma$ is strictly increasing,
the pair $(r,t)$ is uniquely determined
by $\Phi_0(re^{it})$.

In the maximal case, the boundary curve
$\Psi(A+i\R)$ is simple as well.
For sufficiently large $A$, the function
$k(w)=w/\sqrt{\Log w}$ is analytic in a neighborhood
of $\overline{\HH_A}$ and has $\Rea k'(w)>0$ there.
Moreover, $k(\overline{\HH_A})\subset\HH$.
The line-segment argument used above applies to the closed
half-plane and shows that $k$ is injective on
$\overline{\HH_A}$. Since squaring is injective on
$\HH$, the function $\Psi=-k^2$ is injective there
as well; in particular, $\Psi(A+i\R)$ has no
self-intersections.
The two ends of this curve approach the
positive real direction as their moduli
tend to infinity.
By \eqref{eq:polar-expansion},
$\mathcal R(r)$ is strictly increasing
and $\omega(r)$ is strictly decreasing
for sufficiently large $r$.
The image $\Psi(\HH_A)$ is the component
containing the negative real ray.
Outside a sufficiently large disk, the
complementary component containing the
positive real ray is bounded by the
two ends of $\Psi(A+i\R)$.

The polar formula in \eqref{eq:phi-right}
defines a map of
the right half-plane onto this
complementary region near infinity.
Indeed, the image has modulus
$\mathcal R(r)$ and argument between
$-\omega(r)$ and $\omega(r)$.
Since $\mathcal R$ is strictly increasing,
this polar representation is injective.
Thus the left and right images also
have disjoint interiors in the maximal case.

By \eqref{eq:proper-extension}, the
left formula is proper at infinity.
The same property for the right formula
follows from the growth of $\sigma$
and the expressions in \eqref{eq:phi-right}.
Together with the matching boundary values,
the preceding geometric descriptions show
that $\Phi_0$ is a homeomorphism outside
a sufficiently large disk.

We also need local quasiconformality
across the imaginary axis.
Away from the boundary nodes $iLm$,
the slopes of $\sigma$ on the relevant
intervals are bounded above and below
by positive constants.
The differential estimates in
\eqref{eq:local-differential}, together
with the corresponding estimates for
the right-hand formulas, give locally
bounded derivatives and inverse derivatives
on both sides of the imaginary axis.
Since the boundary values agree,
the two formulas define a locally
quasiconformal map across the axis.

It remains to extend $\Phi_0$ across
a compact set. Choose a sufficiently
large radius $R$ that lies strictly
between two consecutive nodes of $\sigma$.
The image of $|z|=R$ under the exterior
map is a simple closed curve.
Its parametrization is piecewise $C^1$,
with nonvanishing one-sided derivatives
and nonzero angles at the junctions.
It follows that the boundary
parametrization is bi-Lipschitz.

The image of the upper semicircle
lies in the upper half-plane, while
the two real endpoints have images
of opposite signs.
Extend the boundary map along
the diameter $[-R,R]$ by an increasing
piecewise affine homeomorphism onto
the segment joining these two images,
chosen to send $0$ to $0$.
This gives a bi-Lipschitz boundary
homeomorphism from the upper half-disk
onto the boundary of a Jordan domain.

Using a fixed bi-Lipschitz identification
of the closed upper half-disk with a
closed disk, we obtain a bi-Lipschitz
embedding of the unit circle.
Lemma~\ref{lem:schoenflies} extends
this embedding to a bi-Lipschitz
homeomorphism of the plane.
Restricting the extension to the
bounded domain gives a bi-Lipschitz
filling of the upper half-disk.
Reflection across the real axis
provides the corresponding extension
on the lower half-disk.

The interior extension agrees with
the exterior map on the bounding
circle and has bounded distortion.
Consequently, the combined map
$\Phi_0$ is a locally quasiconformal
homeomorphism of $\C$.
By construction, $\Phi_0$ commutes
with complex conjugation, fixes
zero and preserves the upper half-plane.

\emph{Step 3. Estimates on a fixed left strip.}
Fix $B>0$ and consider
$-B\le\Rea z\le0$.
For sufficiently large $|w|$ in
$\HH_A$, the comparison functions
in \eqref{eq:comparison-functions}
satisfy
$|\Psi'(w)/\Psi(w)|\le C/|w|$.

Write $z=x+iy$.
The differential estimates
\eqref{eq:local-differential},
\eqref{eq:strip-power} and
\eqref{eq:strip-endpoint} imply
$\|DT(-z)\|\le C_B\sigma(|y|)/|y|$
for sufficiently large $|y|$,
uniformly in $-B\le x\le0$.
Moreover, \eqref{eq:proper-extension}
gives $|A+T(-z)|\asymp\sigma(|y|)$
on this strip.

By \eqref{eq:phi-left} and the
conformality of $\Psi$,
the chain rule therefore gives
\[
 \frac{J_{\Phi_0}(z)}{|\Phi_0(z)|^2}
 \le
 \left|
 \frac{\Psi'(A+T(-z))}
      {\Psi(A+T(-z))}
 \right|^2
 \|DT(-z)\|^2
 \le\frac{C_B}{|y|^2}.
\]
Since $|z|\asymp|y|$ for large
$|z|$ in the fixed strip,
this proves \eqref{eq:comparison-jacobian}
there.

The quasiconformal distortion
of $\Phi_0$ on the strip is
bounded by Lemma~\ref{lem:extension},
because $\Psi$ is conformal.
On the right half-plane,
\eqref{eq:comparison-jacobian}
follows from
\eqref{eq:power-comparison},
\eqref{eq:normal-comparison}
and \eqref{eq:maximal-comparison}.
Together with the bounded distortion
on the compact part of
$\{\Rea z\ge-B\}$, this completes
the proof.
\end{proof}

\subsection{Removing the remaining distortion}

The conformal coordinate $H(z)$ is asymptotic to $-z$
as $\Rea z\to-\infty$. To obtain exact agreement
between the Beltrami coefficients of the comparison
map and the function $\widetilde u$, we modify
$\Phi_0$ on the left half-plane.

Choose $C_0>0$ sufficiently large, and let
$\chi:\R\to[0,1]$ be a smooth function such that
$\chi(x)=1$ for $x\le-2C_0$,
$\chi(x)=0$ for $x\ge-C_0$, and
$|\chi'(x)|\le2/C_0$. Define
\begin{equation}\label{eq:coordinate-correction}
 Q(z)=z-\chi(\Rea z)\Log(1+e^z),\qquad
 \Phi=\Phi_0\circ Q.
\end{equation}
The logarithm is defined by its convergent power
series wherever $\chi$ is nonzero, and the
correction vanishes for $\Rea z\ge-C_0$.
Since $|e^z|\le e^{-C_0}$ on the support
of $\chi$, the correction in
\eqref{eq:coordinate-correction} has Lipschitz
constant less than $1/2$ when $C_0$ is
sufficiently large. Thus $Q$ is a
bi-Lipschitz homeomorphism of $\C$.
It commutes with complex conjugation and fixes zero.

For $\Rea z\le-2C_0$, we have
$Q(z)=z-\Log(1+e^z)=-H(z)$
by \eqref{eq:left-coordinate}.
Consequently, \eqref{eq:phi-left} gives,
for sufficiently large $|z|$ in this region,
\begin{equation}\label{eq:exact-agreement}
 \Phi(z)=\Psi(A+T(H(z))),\qquad
 \mu_\Phi(z)=\mu(z).
\end{equation}
The equality of the Beltrami coefficients
follows because $\Psi$ and $H$ are conformal,
while the coefficient of $T\circ H$ is
precisely $\mu$ by \eqref{eq:physical-mu}.

Since $Q$ is bi-Lipschitz and has bounded
displacement, the Jacobian estimate
\eqref{eq:comparison-jacobian} also holds
for $\Phi=\Phi_0\circ Q$ on every fixed strip
$-B\le\Rea z\le0$. On the right half-plane,
$Q$ is the identity, so $\Phi=\Phi_0$ there.

\begin{prop}\label{prop:straightening}
There exists an orientation-preserving,
locally quasiconformal homeomorphism
$g:\C\to\C$ with Beltrami coefficient
$\mu$ defined by \eqref{eq:physical-mu}.
The map $g$ commutes with complex conjugation,
fixes zero and preserves the upper half-plane.
Moreover,
\begin{equation}\label{eq:g-asymptotic}
 g(z)=\Phi(z)+o(|\Phi(z)|)\qquad(z\to\infty)
\end{equation}
uniformly in direction.
\end{prop}

\begin{proof}
We seek $g$ in the form $G\circ\Phi$,
where $G$ is a global quasiconformal
homeomorphism. To determine the Beltrami
coefficient of $G$, introduce the
coordinate $W=\Phi(z)$ and define
\begin{equation}\label{eq:relative-coefficient}
 \nu(\Phi(z))=
 \frac{\Phi_z(z)}{\overline{\Phi_z(z)}}\,
 \frac{\mu(z)-\mu_\Phi(z)}
      {1-\overline{\mu_\Phi(z)}\mu(z)}.
\end{equation}
The composition formula for Beltrami
coefficients shows that if $G$ has
coefficient $\nu$, then $G\circ\Phi$
has coefficient $\mu$.

We first prove that $\|\nu\|_\infty<1$.
By \eqref{eq:exact-agreement}, the
numerator in \eqref{eq:relative-coefficient}
vanishes when $\Rea z\le-2C_0$ and
$|z|$ is sufficiently large.
Thus $\nu$ vanishes on the image
under $\Phi$ of this region.

On the half-plane $\Rea z\ge-2C_0$,
the coefficients $\mu$ and $\mu_\Phi$
are uniformly bounded away from modulus
one. For $\mu$, this follows from
Lemma~\ref{lem:extension} and
\eqref{eq:physical-mu}. For $\mu_\Phi$,
it follows from Lemma~\ref{lem:comparison}
and the bi-Lipschitz property of $Q$.
Both coefficients are also bounded
away from modulus one on the remaining
compact region, since the corresponding
maps are locally quasiconformal.
The disk automorphism appearing in
\eqref{eq:relative-coefficient} therefore
gives $\|\nu\|_\infty<1$.

We next verify the integrability condition
\begin{equation}\label{eq:relative-integral}
 \int_{|W|>R}\frac{|\nu(W)|}{|W|^2}
 \dd A(W)<\infty
\end{equation}
for sufficiently large $R$.

Changing variables by $W=\Phi(z)$,
we have to estimate
$|\nu(\Phi(z))|J_\Phi(z)/|\Phi(z)|^2$.
By \eqref{eq:exact-agreement}, this
quantity vanishes for sufficiently
large $|z|$ when $\Rea z\le-2C_0$.

On the strip $-2C_0\le\Rea z\le0$,
the coefficient $\nu$ is uniformly
bounded, and the Jacobian estimate
\eqref{eq:comparison-jacobian},
which also holds for $\Phi$, gives
\[
 \frac{|\nu(\Phi(z))|J_\Phi(z)}
      {|\Phi(z)|^2}
 \le\frac{C}{|z|^2}.
\]
The right-hand side is integrable
outside a disk on this fixed strip.

For $\Rea z\ge0$, we have
$\Phi=\Phi_0$. Since both $\mu$
and $\mu_{\Phi_0}$ are uniformly
bounded away from modulus one,
\eqref{eq:relative-coefficient}
and \eqref{eq:comparison-jacobian}
give
\begin{equation}\label{eq:relative-majorant}
 C|\mu(z)|\,|z|^{-2}
   +C|\mu_{\Phi_0}(z)|\,
                 J_{\Phi_0}(z)/|\Phi_0(z)|^2
\end{equation}
as an upper bound for the transformed
integrand. The first term in
\eqref{eq:relative-majorant} is integrable
by \eqref{eq:physical-integrability},
and the second is integrable by
\eqref{eq:comparison-integral}.
This proves \eqref{eq:relative-integral}.

Since $\|\nu\|_\infty<1$,
Lemma~\ref{lem:mrmt} provides a unique
global quasiconformal homeomorphism
$G:\C\to\C$ with Beltrami coefficient
$\nu$, normalized by $G(0)=0$,
$G(1)=1$ and $G(\infty)=\infty$.
The integrability condition
\eqref{eq:relative-integral} allows
us to apply Lemma~\ref{lem:twb},
which gives $G(W)=cW+o(|W|)$
as $W\to\infty$, for some $c\ne0$.

The symmetry of $\mu$ and $\Phi$
under complex conjugation implies
that $\nu$ has the same symmetry.
By uniqueness of the normalized
solution, $G$ commutes with complex
conjugation. Hence $G$ maps the
real axis onto itself. Since $G$
is orientation-preserving and fixes
$0$ and $1$, its restriction to
the real axis is increasing.
It follows that $c>0$.

Finally, define $g=c^{-1}G\circ\Phi$.
The composition formula and
\eqref{eq:relative-coefficient}
show that $g$ has Beltrami
coefficient $\mu$. Since
$\Phi(z)\to\infty$ as $z\to\infty$,
the asymptotic formula for $G$
gives \eqref{eq:g-asymptotic}.
The map $g$ commutes with complex
conjugation and fixes zero.
Moreover, both $G$ and $\Phi$
preserve the upper half-plane,
and multiplication by $c^{-1}>0$
preserves it as well.
This completes the proof.
\end{proof}

\section{The potential and its asymptotics}\label{sec:completion}

Let $g$ be the homeomorphism constructed in
Proposition~\ref{prop:straightening}, and define
$v=\widetilde u\circ g^{-1}$.
By Lemma~\ref{lem:gluing}, the function $v$ is subharmonic
and has no critical points outside its logarithmic poles.
The poles in the upper half-plane are
\begin{equation}\label{eq:sites}
 z_m=g((2m-1)\pi i),\qquad m\ge1.
\end{equation}
Since $g$ commutes with complex conjugation, the remaining
poles are $\overline{z_m}$, $m\ge1$. All the poles are
distinct and nonzero because $g$ is a homeomorphism
fixing zero. By \eqref{eq:slopes} and
Lemma~\ref{lem:gluing}, the logarithmic coefficient
at both $z_m$ and $\overline{z_m}$ is $b_m$.

We first describe the growth of the positive part of $v$.
For the corresponding choice of $A$ and $\Psi$ in
\eqref{eq:comparison-functions}, put
$\mathcal D=\Psi(\HH_A)$.
Define a nonnegative function $P$ on $\C$ by
\[
 P(Z)=
 \begin{cases}
  \Rea\Psi^{-1}(Z)-A,&Z\in\mathcal D,\\
  0,&Z\notin\mathcal D.
 \end{cases}
\]
Here $\Psi^{-1}$ denotes the conformal inverse of
$\Psi:\HH_A\to\mathcal D$. The function $P$ is
continuous on $\C$ and vanishes on the boundary
of $\mathcal D$.

\begin{lem}\label{lem:model-growth}
In the power case, $P=h_\rho$, where $h_\rho$ is
defined by \eqref{eq:angular-profile}.

In both the normal and maximal cases,
\begin{equation}\label{eq:model-angular}
 \frac{P(re^{it})}{S(r)}\longrightarrow|\sin(t/2)|
 \quad\text{uniformly for }-\pi\le t\le\pi,
\end{equation}
as $r\to\infty$, with the corresponding growth
scale $S$ from \eqref{eq:endpoint-models}.

Moreover, in all three cases,
\begin{equation}\label{eq:positive-pullback}
 \widetilde u^+(z)=P(\Phi(z))+O(1)\qquad(|z|\to\infty)
\end{equation}
uniformly in direction.
\end{lem}

\begin{proof}
In the power case, $\Psi(w)=-w^q$ with
$q=1/\rho$. Inverting this map on $\mathcal D$
and taking the real part gives
$P(Z)=h_\rho(Z)$, as defined in
\eqref{eq:angular-profile}.

In the normal case, $\Psi(w)=-w^2$ and
$A=\pi/2$. Taking the square root with
nonnegative real part, we obtain
$P(Z)=(\Rea\sqrt{-Z}-A)^+$.
Consequently,
$P(re^{it})/\sqrt r=|\sin(t/2)|+O(r^{-1/2})$,
uniformly for $t\in[-\pi,\pi]$.
This proves \eqref{eq:model-angular} in
the normal case.

For the maximal case, let $Z\in\mathcal D$
and write $Z=-w^2/\Log w$, where
$w=\Psi^{-1}(Z)\in\HH_A$.
Uniformly as $|w|\to\infty$ in $\HH_A$,
\begin{equation}\label{eq:inverse-asymptotic}
 |w|\sim\sqrt{|Z|\log|Z|/2},\qquad
 \arg w=\tfrac12\arg(-Z)+O(1/\log|Z|).
\end{equation}
Indeed, taking absolute values in
$Z=-w^2/\Log w$ gives
$|Z|=|w|^2/|\Log w|$ and hence
$\log|Z|=2\log|w|-\log\log|w|+o(1)$.
Moreover, $|\arg\Log w|=O(1/\log|w|)$,
which gives the angular estimate in
\eqref{eq:inverse-asymptotic}.
Here $\arg(-Z)$ is chosen in $(-\pi,\pi)$
on $\mathcal D$.

Since $P(Z)=\Rea w-A$ for $Z\in\mathcal D$,
taking real parts in \eqref{eq:inverse-asymptotic}
gives \eqref{eq:model-angular} on $\mathcal D$,
with $S(r)=\sqrt{r\log r/2}$.
For $|Z|=r$ outside $\mathcal D$, the
boundary-angle estimate in
\eqref{eq:polar-expansion} shows that
$|\arg Z|=O(1/\log r)$ for sufficiently
large $r$. Thus $|\sin(t/2)|$ tends
uniformly to zero on this part of the
circle, while $P(Z)=0$ by definition.
This proves \eqref{eq:model-angular}
uniformly on the whole circle.

It remains to establish \eqref{eq:positive-pullback}.
For $\Rea z\le-2C_0$ and sufficiently large
$|z|$, the exact identity
\eqref{eq:exact-agreement} gives
$\Phi(z)=\Psi(A+T(H(z)))$.
Since $T(H(z))\in\HH$, the definition
of $P$ yields
$P(\Phi(z))=\Rea T(H(z))
=\widetilde u^+(z)$.

On the half-plane $\Rea z\ge-2C_0$,
the positive part $\widetilde u^+$ is
uniformly bounded. Indeed, on $D$,
\eqref{eq:coordinate-bounds} bounds
$\Rea H(z)$, and Lemma~\ref{lem:extension}
gives $\Rea T(H(z))\le C\Rea H(z)$.
Outside $D$, we have $\widetilde u\le0$.

We also need to bound $P(\Phi(z))$ on
this half-plane. On the right,
$Q$ is the identity and the image of
\eqref{eq:phi-right} lies outside
$\mathcal D$, so $P(\Phi(z))=0$
for sufficiently large $|z|$.
On a fixed strip to the left of
the imaginary axis, \eqref{eq:phi-left}
gives $P(\Phi_0(z))=\Rea T(-z)
\le C|\Rea z|$.
The bi-Lipschitz map $Q$ in
\eqref{eq:coordinate-correction} has
bounded displacement, and the same
bound therefore holds for $P(\Phi(z))$
on the strip.

Thus $\widetilde u^+(z)$ and
$P(\Phi(z))$ differ by a bounded
quantity outside a sufficiently large
disk, proving \eqref{eq:positive-pullback}.
\end{proof}

We now prove Theorems~\ref{thm:main} and
\ref{thm:endpoint}, together with the
growth assertions in Corollary~\ref{cor:growth}.

\begin{proof}
\emph{Step 1. Pole locations and summability.}
Recall that the upper half-plane poles
are given by \eqref{eq:sites}.
Since $Q$ is the identity on the imaginary
axis, \eqref{eq:phi-left} gives
$\Phi((2m-1)\pi i)
=\Psi(A-i\sigma((2m-1)\pi))$
for sufficiently large $m$.
Combining this identity with
\eqref{eq:g-asymptotic}, we obtain
\begin{equation}\label{eq:pole-model}
 z_m\sim\Psi(A-i\sigma((2m-1)\pi)).
\end{equation}

In the power case, \eqref{eq:power-boundary}
gives
$\sigma((2m-1)\pi)\sim(L/\beta)m^\beta$.
Since $\Psi(w)=-w^q$ and $p=\beta q=2-q$,
\eqref{eq:pole-model} becomes
\begin{equation}\label{eq:pole-asymptotic}
 z_m\sim\kappa m^p e^{i\pi p/2},
 \qquad \kappa=(L/\beta)^q .
\end{equation}

The two asymptotic rays of poles are illustrated in
Figure~\ref{fig:pole-angular-growth}(a) for $\rho=3/4$.
By \eqref{eq:mass-profile},
$b_m=O(m^{\beta-1})$.
Hence $b_m/|z_m|=O(m^{\beta-1-p})$.
Since $p>\beta$, the series
$\sum_{m\ge1}b_m/|z_m|$ converges.

In the normal and maximal cases,
\eqref{eq:endpoint-boundary} and
\eqref{eq:pole-model} give the pole
asymptotics in \eqref{eq:endpoint-poles}.
Using the masses in \eqref{eq:endpoint-masses},
we obtain the respective bounds
\[
 \frac1{m(\log m)^2}\quad\text{and}\quad
 \frac1{m\exp\sqrt{\log m}},
\]
up to constant factors, for $b_m/|z_m|$.
Both series converge. For the second,
this follows from the integral test with
the substitution $t=\sqrt{\log m}$.
Thus the summability condition
\eqref{eq:sum} holds in all three cases.

\emph{Step 2. Identification of the potential
and its positive part.}
By Proposition~\ref{prop:straightening},
$g(z)=\Phi(z)+o(|\Phi(z)|)$ uniformly
as $z\to\infty$.
Consequently, if $Z=g(z)$, then
$\Phi(z)=Z+o(|Z|)$ uniformly in direction.

Recall that $v=\widetilde u\circ g^{-1}$.
In the power case, Lemma~\ref{lem:model-growth}
gives $P=h_\rho$; the homogeneity and
continuity of $h_\rho$ then allow us to
transfer \eqref{eq:positive-pullback}
from the $\Phi$-coordinate to the
$Z$-coordinate. In the endpoint cases,
the same transfer follows from
\eqref{eq:model-angular} and the
regular variation of $S$.
We obtain
\begin{equation}\label{eq:straightened-profile}
 \begin{aligned}
 v^+(Z)&=h_\rho(Z)+o(|Z|^\rho)
                              &&\text{in the power case},\\
 v^+(re^{it})&=S(r)|\sin(t/2)|+o(S(r))
                              &&\text{in either endpoint case},
 \end{aligned}
\end{equation}
where the error estimates are uniform
in the angular variable.

For completeness, both endpoint growth
scales satisfy $S(cr)/S(r)\to\sqrt c$
as $r\to\infty$, uniformly for $c$
in compact subintervals of $(0,\infty)$.
Thus replacing $r$ by $r(1+o(1))$
changes $S(r)$ by a factor $1+o(1)$.
Together with the continuity of
$|\sin(t/2)|$, this justifies the
endpoint assertions in
\eqref{eq:straightened-profile}.

Since $\rho<1$ and $S(r)=o(r)$ in both endpoint cases,
\eqref{eq:straightened-profile} implies
\[
 \sup_{|Z|=r}v^+(Z)=o(r).
\]
The poles and their logarithmic
coefficients have already been
identified in Step~1 and
Lemma~\ref{lem:gluing}.
We may therefore apply
Lemma~\ref{lem:harmonic-remainder}
to obtain
\begin{equation}\label{eq:exact-potential}
 u(Z)=v(Z)-v(0)=\sum_{m\ge1}b_m
       \bigl(\log|1-Z/z_m|+\log|1-Z/\overline{z_m}|\bigr).
\end{equation}

The summability established in Step~1
justifies termwise differentiation
away from the poles. Thus the
conjugate gradient of $u$ is exactly
\eqref{eq:main-series}.
By Lemma~\ref{lem:gluing}, this gradient
has no zeros.

Since $u=v-v(0)$, their positive parts
differ by at most $|v(0)|$.
Hence \eqref{eq:straightened-profile}
also gives the corresponding asymptotic
formulas for $u^+$.

In the power case, the value of
$\kappa$ in \eqref{eq:pole-asymptotic}
gives $D_\rho=2\kappa^{-\rho}/\beta=1/\pi$.
Thus \eqref{eq:straightened-profile}
proves \eqref{eq:angular-limit}
with this normalization.
A dilation of the $Z$-coordinate
changes $\kappa$ to any prescribed
positive value without changing
the masses, and gives the factor
$\pi D_\rho$ in \eqref{eq:angular-limit}.
In the endpoint cases,
\eqref{eq:straightened-profile}
proves \eqref{eq:endpoint-angular}.

Figure~\ref{fig:pole-angular-growth}(b) compares the normalized angular profiles.
Together with the pole asymptotics from Step~1, it shows
that the asymptotic rays of poles bound the sector where the leading-order
growth of $u^+$ vanishes.

\begin{figure}[htbp]
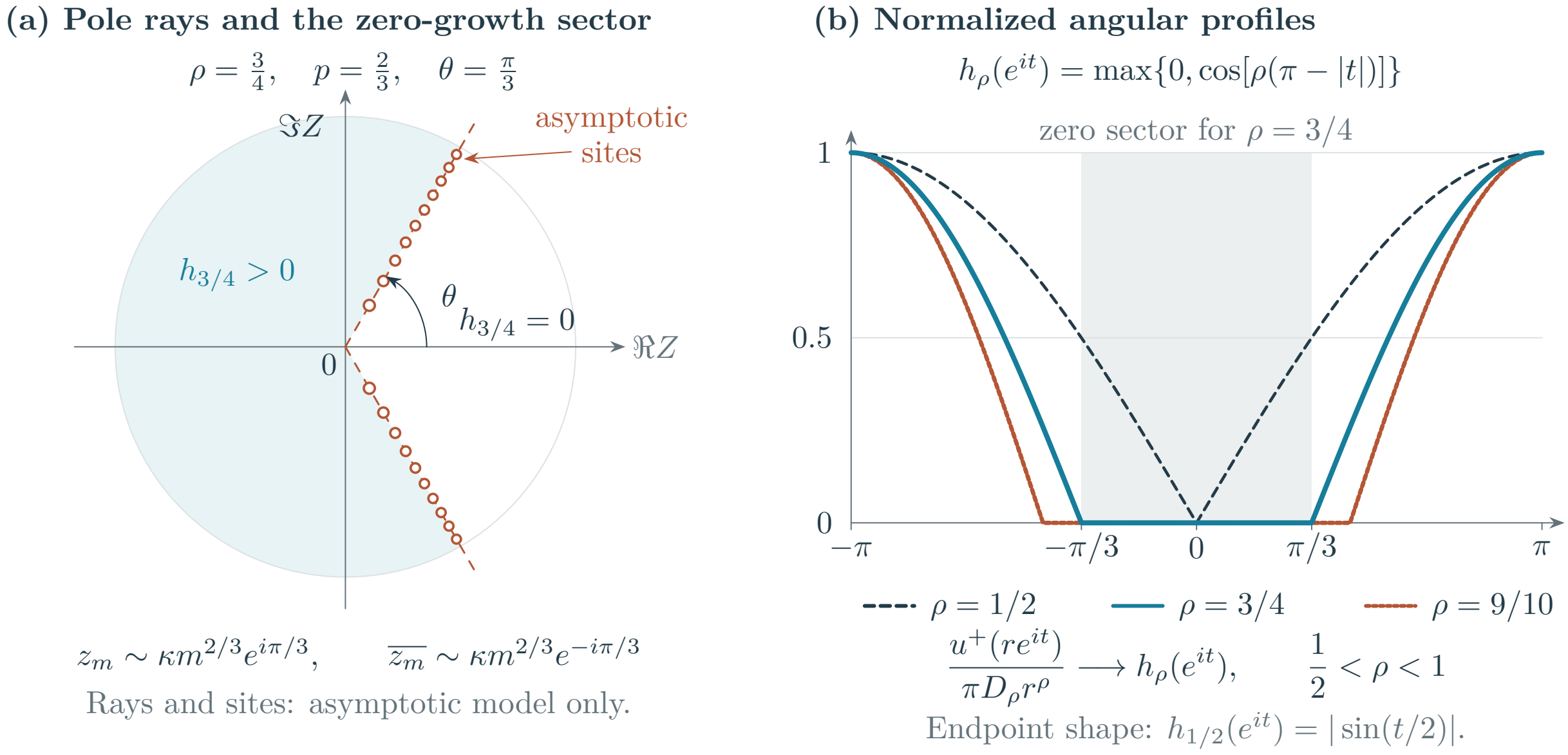

\centering
\begingroup
\resizebox{\linewidth}{!}{%
%
}
\endgroup
\caption{Limiting pole geometry and angular growth.
(a) For $\rho=3/4$, the poles approach the rays
$\arg z=\pm\pi/3$; the displayed sites represent the asymptotic model.
(b) The normalized angular profiles $h_\rho(e^{it})$ for
$\rho=1/2,3/4,9/10$. For $\rho>1/2$, the normalized positive part
converges as in \eqref{eq:angular-limit}. The curve for $\rho=1/2$
is the endpoint shape $|\sin(t/2)|$, with the normalization $S(r)$
from \eqref{eq:endpoint-angular}. The shaded interval and sector
record vanishing leading-order growth for $\rho=3/4$.}
\label{fig:pole-angular-growth}
\end{figure}

\emph{Step 3. Mass, maxima and circular means.}
We first consider the power case.
By \eqref{eq:sigma} and
\eqref{eq:power-boundary}, the
cumulative masses satisfy
$\sum_{m\le k}b_m\sim k^\beta/\beta$.
On the other hand,
\eqref{eq:pole-asymptotic} gives
$|z_m|\sim\kappa m^p$.

Fix $0<\varepsilon<1$.
For sufficiently large $m$,
$(1-\varepsilon)\kappa m^p
\le|z_m|\le
(1+\varepsilon)\kappa m^p$.
Thus the indices for which
$|z_m|\le r$ can be bounded
using the cutoffs
$(r/((1+\varepsilon)\kappa))^{1/p}$
and
$(r/((1-\varepsilon)\kappa))^{1/p}$.
Applying the cumulative mass
asymptotic and then letting
$\varepsilon\downarrow0$ gives
\begin{equation}\label{eq:power-mass-growth}
 M(r)\sim D_\rho r^\rho,\qquad
 D_\rho=2\kappa^{-\rho}/\beta.
\end{equation}
Here the factor $2$ accounts
for the conjugate pair of poles
carrying each mass $b_m$.

For the endpoint cases, introduce
the model radii
$R_m=|\Psi(A-i\sigma((2m-1)\pi))|$.
The sequence $(R_m)$ is eventually
increasing, and
\eqref{eq:endpoint-boundary} gives
\begin{equation}\label{eq:endpoint-cumulative}
 2\sum_{j\le m}b_j=\frac{2h(Lm)}{L}
       \sim\frac{S(R_m)}{\pi},\qquad R_{m+1}/R_m\to1.
\end{equation}
In the normal case, these relations
follow directly from
$h(t)=\log t$ and $\Psi(w)=-w^2$.

For the maximal case, put
$\ell=\log h(Lm)$.
The definition of $\Psi$ gives
$R_m\sim e^{2\ell}/\ell$
and $\log R_m\sim2\ell$.
Substituting these asymptotics into
$S(r)=\sqrt{r\log r/2}$
proves \eqref{eq:endpoint-cumulative}.

Since $R_m$ is eventually increasing
and $R_{m+1}/R_m\to1$, the cumulative
mass obtained by counting the
model radii $R_m$ is asymptotic
to $S(r)/\pi$.
More precisely,
$2\sum_{R_m\le r}b_m\sim S(r)/\pi$,
by \eqref{eq:endpoint-cumulative}
and the regular variation of $S$.

By \eqref{eq:pole-model}, the actual
pole moduli satisfy $|z_m|\sim R_m$.
For every fixed $0<\varepsilon<1$,
we therefore have
$(1-\varepsilon)R_m
\le|z_m|\le(1+\varepsilon)R_m$
for all sufficiently large $m$.
Comparing the cumulative masses
at these two radii and letting
$\varepsilon\downarrow0$ gives
$M(r)\sim S(r)/\pi$, as asserted
in \eqref{eq:endpoint-means}.

This argument does not require
the sequence $(|z_m|)$ itself
to be monotone.

The uniform limits for $u^+$, and the subharmonic maximum principle,
give $B(r,u)\sim\pi D_\rho r^\rho$ in the power case and
$B(r,u)\sim S(r)$ in either endpoint case.

We next determine the circular means. Write
\[
 \mathcal M_+(r,u)=\frac1{2\pi}\int_{-\pi}^{\pi}u^+(re^{it})\dd t.
\]
The mean of $h_\rho$ on the unit
circle is $1/(\pi\rho)$, while
the mean of $|\sin(t/2)|$ is $2/\pi$.
Integrating the uniform limits in
\eqref{eq:angular-limit} and
\eqref{eq:endpoint-angular} therefore
gives $\mathcal M_+(r,u)\sim D_\rho r^\rho/\rho$
in the power case and
$\mathcal M_+(r,u)\sim2S(r)/\pi$
in the endpoint cases.

To obtain the corresponding
$L^1$ limits for $u$, rather
than only for $u^+$, we use
Jensen's formula:
\begin{equation}\label{eq:jensen}
 \frac1{2\pi}\int_{-\pi}^{\pi}u(re^{it})\dd t
       =\int_0^r\frac{M(s)}s\dd s .
\end{equation}
In the power case,
\eqref{eq:power-mass-growth} gives
$\int_0^r M(s)\dd s/s
\sim D_\rho r^\rho/\rho$.

In the endpoint cases, the mass
asymptotic obtained from
\eqref{eq:endpoint-cumulative}
gives
$\int_0^r M(s)\dd s/s\sim2S(r)/\pi$.
To justify this last equivalence,
substitute $s=rt$.
For $rt$ above a fixed large
constant, both endpoint scales
satisfy
$S(rt)/S(r)\le Ct^{1/2}$
when $0<t\le1$.
Dominated convergence then
gives the stated asymptotic
for the integral.

In every case, the mean of $u$
in \eqref{eq:jensen} has the same
leading term as the mean of $u^+$.
Since $u=u^+-u^-$, where
$u^-=\max\{-u,0\}$, the mean
of $u^-$ is therefore of smaller
order than $r^\rho$ in the power
case and smaller than $S(r)$
in the endpoint cases.
Logarithmic singularities lying
on a circle are integrable in
the angular variable.
Consequently, the uniform
limits for $u^+$ imply the
$L^1[-\pi,\pi]$ limits for $u$
stated in Corollary~\ref{cor:growth}
and Theorem~\ref{thm:endpoint}.

\emph{Step 4. The meromorphic order
in the power case.}
All poles of the function
$f$ in \eqref{eq:main-series}
are simple.
By \eqref{eq:pole-asymptotic},
their ordinary counting function
satisfies
$n_p(r)\sim2\kappa^{-1/p}r^{1/p}$.
Lemma~\ref{lem:proximity} gives
$m(r,f)=o(1)$, so the Nevanlinna
characteristic satisfies
\begin{equation}\label{eq:meromorphic-characteristic}
 T(r,f)\sim2p\kappa^{-1/p}r^{1/p}.
\end{equation}
Hence both the order and lower
order of $f$ are $1/p$, as
asserted in Corollary~\ref{cor:growth}.
\end{proof}

The asymptotic pole rays in
Theorem~\ref{thm:main} also show
that an angular restriction of
Langley and Rossi~\cite[Theorem~1.6]{LR04}
is sharp.
Their theorem implies that a positive Cauchy sum $f$
satisfying \eqref{eq:sum},
with positive residues of infinite
total mass and finite lower order
$\lambda>1$, has infinitely many zeros
provided that
$\limsup_j|\arg\zeta_j|<\pi/(2\lambda)$.
We apply the case of parameter value $1$ in
\cite[Theorem~1.6]{LR04}; finitely many
exceptional poles may be absorbed
into the rational term allowed there.

For the functions $f$ constructed in
the power case, \eqref{eq:meromorphic-characteristic}
gives $\lambda=1/p$.
Moreover, \eqref{eq:pole-asymptotic}
shows that the poles approach the rays
$\arg z=\pm\pi/(2\lambda)$.
Thus our zero-free examples attain
the boundary angle in the theorem
of Langley and Rossi.

\section{Finite total mass and analytic capacity}
\label{sec:finite-mass}

For a compact set $E\subset\C$, let $\gamma(E)$ denote
its analytic capacity, defined as the supremum of
$|\lim_{z\to\infty}zF(z)|$ over all functions $F$
holomorphic on $\C\setminus E$ satisfying
$|F|\le1$ and $F(\infty)=0$.
For an arbitrary set $E\subset\C$, define
$\gamma(E)$ as the supremum of $\gamma(K)$
over compact subsets $K\subset E$.

For a finite complex measure $\nu$, let $\|\nu\|$
denote its total variation. The maximal truncated
Cauchy transform of $\nu$ is defined by
\[
 C_*\nu(z)=\sup_{\varepsilon>0}
 \left|\int_{|z-t|>\varepsilon}
 \frac{\dd\nu(t)}{z-t}\right|.
\]
We use the weak capacity estimate for the Cauchy
transform and the countable semiadditivity of
analytic capacity due to Tolsa~\cite{Tol02,Tol03},
in the formulations given by Mateu, Tolsa and
Verdera~\cite[Theorems~1.1--1.2 and Section~2]{MTV03}.
We also use the diameter estimate for continua
stated by Totik~\cite[Section~2, formula~(5),
p.~199]{Tot22}.

\begin{lem}[{\cite[Theorems~1.1--1.2 and Section~2]{MTV03},
\cite[formula~(5), p.~199]{Tot22}}]
\label{lem:capacity-facts}
There exists an absolute constant $C>0$ such that,
for every finite complex measure $\nu$ and every
$\lambda>0$,
\begin{equation}\label{eq:weak-capacity}
 \gamma\{z:C_*\nu(z)>\lambda\}
 \le C\|\nu\|/\lambda.
\end{equation}
Moreover, for every sequence $(E_j)_{j\ge1}$
of Borel subsets of $\C$,
$\gamma(\bigcup_{j\ge1}E_j)
\le C\sum_{j\ge1}\gamma(E_j)$.
Every countable set has analytic capacity zero,
and every compact connected set $K\subset\C$
satisfies $\operatorname{diam}K\le4\gamma(K)$.
\end{lem}

To compare the zeros and poles of $f$,
we need a Jordan curve surrounding a given disk
while avoiding a compact exceptional set.
The following lemma provides such a curve when
the exceptional set has sufficiently small
analytic capacity.

\begin{lem}\label{lem:annular-circuit}
Let $0<r<R$, and let $E$ be a compact subset
of the annulus $\{z:r\le|z|\le R\}$.
If $\gamma(E)<(R-r)/4$, then there exists
a polygonal Jordan curve
$\Gamma\subset\{z:r<|z|<R\}\setminus E$
whose bounded complementary component contains
the closed disk $\{z:|z|\le r\}$.
\end{lem}

\begin{proof}
Choose $d$ such that
$4\gamma(E)<d<R-r$.
By Lemma~\ref{lem:capacity-facts},
every connected component $K$ of $E$
satisfies
$\operatorname{diam}K\le4\gamma(K)
\le4\gamma(E)<d$.

Consider a square grid of mesh size $\delta$,
and let $E_\delta$ be the union of all closed
grid squares meeting $E$.
We claim that, for sufficiently small $\delta$,
every connected component of $E_\delta$
has diameter less than $d$.

Suppose otherwise. Then there are mesh sizes
$\delta_n\to0$ and connected components
$K_n$ of $E_{\delta_n}$ with
$\operatorname{diam}K_n\ge d$.
Since $E$ is compact and every point of
$E_{\delta_n}$ lies within distance
$\sqrt2\delta_n$ of $E$, the sets $K_n$
lie in a fixed compact disk.
By compactness of the hyperspace of nonempty compact
subsets of a fixed compact disk in the Hausdorff metric,
a subsequence of $(K_n)$ converges to a nonempty compact
set $K$. This limit is connected: if two disjoint open
sets separated $K$, they would also separate $K_n$
for all sufficiently large $n$, a contradiction.
Since every point of $K_n$ lies within
$\sqrt2\delta_n$ of $E$, we have $K\subset E$;
moreover, $\operatorname{diam}K\ge d$.
This contradicts the diameter bound for connected
subsets of $E$.
The claim follows.

Fix a grid for which every component of
$E_\delta$ has diameter less than $d$.
Let $K_0$ be the union of the closed disk
$\{|z|\le r\}$ and all components of
$E_\delta$ meeting that disk.
Let $K_1$ be the union of the remaining
components of $E_\delta$ and the closed
set $\{|z|\ge R\}$.

Since $d<R-r$, no component of $E_\delta$
can meet both $\{|z|\le r\}$ and
$\{|z|\ge R\}$.
Thus $K_0$ and $K_1$ are disjoint.
Moreover, $K_0$ is compact and $K_1$
is closed, so $\dist(K_0,K_1)>0$.

Choose a triangulation of the plane whose
mesh is less than one quarter of
$\dist(K_0,K_1)$.
Assign value $1$ to every vertex within
one mesh length of $K_0$, and value $0$
to all other vertices.
Extend these values affinely over each
triangle. The resulting continuous
piecewise affine function has compact
support, equals $1$ on $K_0$, and
vanishes on $K_1$.

Since no vertex has value $1/2$,
the level set at $1/2$ consists of
finitely many disjoint polygonal
Jordan curves. Each triangle meeting
this level set is intersected in a
single line segment.

Consider the connected component of
the set where the function exceeds
$1/2$ that contains $\{|z|\le r\}$.
Its outer boundary is a polygonal
Jordan curve $\Gamma$.
The curve lies between $K_0$ and
$K_1$, so it avoids $E_\delta$
and is contained in
$\{z:r<|z|<R\}$.
Since $E\subset E_\delta$,
the curve also avoids $E$.
By construction, the bounded interior
of $\Gamma$ contains the closed disk
$\{|z|\le r\}$.
\end{proof}

The constant $1/4$ in
Lemma~\ref{lem:annular-circuit} is sharp.
Indeed, the radial segment $E=[r,R]$
meets every Jordan curve in the annulus
whose bounded interior contains
$\{|z|\le r\}$.

To calculate its analytic capacity, put
$c=(r+R)/2$ and $h=(R-r)/2$.
Choose the branch of
$\sqrt{(z-c)^2-h^2}$ that is
asymptotic to $z-c$ at infinity.
The function
$z\mapsto h/(z-c+\sqrt{(z-c)^2-h^2})$
maps the complement of $E$ in the
Riemann sphere conformally onto the
unit disk. Its Laurent expansion at
infinity begins with $h/(2z)$.
Schwarz's lemma, applied after
composition with the inverse map,
shows that no admissible function
has a larger coefficient at infinity.
Consequently,
$\gamma([r,R])=h/2=(R-r)/4$,
as claimed.

We now prove Theorem~\ref{thm:finite-mass}
by comparing $f$ with $A/z$,
where $A=\sum_{j\ge1}a_j\ne0$.

\begin{proof}
Fix $\eta>0$ and put
$B=\sum_{j\ge1}|a_j|$.
For sufficiently large $r$, split $f$
into
$f=g_r+h_r$, where
$g_r(z)=\sum_{|\zeta_j|\le\sqrt r}
a_j/(z-\zeta_j)$ and
$h_r(z)=\sum_{|\zeta_j|>\sqrt r}
a_j/(z-\zeta_j)$.
Write
$\tau_r=\sum_{|\zeta_j|>\sqrt r}|a_j|$.
Since the residues are absolutely
summable, $\tau_r\to0$ as $r\to\infty$.

Consider the closed annulus
$A_r=\{z:r\le|z|\le(1+\eta)r\}$.
For $z\in A_r$, subtracting $A/z$
from $g_r(z)$ gives
\begin{equation}\label{eq:head-approximation}
 \left|g_r(z)-\frac A z\right|
 \le \frac{B\sqrt r}{r(r-\sqrt r)}
       +\frac{\tau_r}r
 =o(r^{-1}).
\end{equation}
Indeed, for $|\zeta_j|\le\sqrt r$,
we have
$|1/(z-\zeta_j)-1/z|
\le\sqrt r/[r(r-\sqrt r)]$.
The remaining difference comes from
replacing the sum of the retained
residues by the full sum $A$,
and its modulus is at most $\tau_r/r$.
Thus the estimate in
\eqref{eq:head-approximation}
holds uniformly on $A_r$.

We next control the tail $h_r$.
Put $\lambda_r=|A|/[4(1+\eta)r]$,
and define
$E_r=\{z\in A_r:|h_r(z)|\ge\lambda_r\}$,
including the poles of $h_r$.
Since $h_r$ is meromorphic and has
only finitely many poles in $A_r$,
the set $E_r$ is compact.

Consider the finite complex measure
$\nu_r=\sum_{|\zeta_j|>\sqrt r}
a_j\delta_{\zeta_j}$, whose total
variation is $\tau_r$.
Away from the poles, the Cauchy
transform of $\nu_r$ equals $h_r$,
and $|h_r(z)|\le C_*\nu_r(z)$.
Hence the part of $E_r$ outside
the poles is contained in
$\{z:C_*\nu_r(z)>\lambda_r/2\}$.
The poles form a countable set
of analytic capacity zero.
Applying \eqref{eq:weak-capacity}
and the semiadditivity estimate
in Lemma~\ref{lem:capacity-facts},
we obtain
\begin{equation}\label{eq:small-exceptional-capacity}
 \gamma(E_r)
 \le C\tau_r/\lambda_r=o(r).
\end{equation}

Since the width of $A_r$ is
$\eta r$, \eqref{eq:small-exceptional-capacity}
gives $\gamma(E_r)<\eta r/4$
for sufficiently large $r$.
Lemma~\ref{lem:annular-circuit}
therefore provides a positively
oriented polygonal Jordan curve
$\Gamma_r\subset A_r\setminus E_r$
whose interior contains
$\{|z|\le r\}$.

The choice of $\Gamma_r$ is illustrated in
Figure~\ref{fig:finite-mass-circuit}(a).

On $\Gamma_r$, the definition
of $E_r$ gives $|h_r(z)|<\lambda_r$.
By \eqref{eq:head-approximation},
we also have
$|g_r(z)-A/z|<\lambda_r$
for sufficiently large $r$.
Consequently,
$|f(z)-A/z|<2\lambda_r
\le|A/z|/2$ on $\Gamma_r$.

The homotopy
$A/z+t(f(z)-A/z)$, $0\le t\le1$,
therefore has no zeros on
$\Gamma_r$.
Since $\Gamma_r$ surrounds the
origin once in the positive
direction, the winding number
of $A/z$ along $\Gamma_r$ is
$-1$. The same is therefore
true of $f$.

The homotopy and the resulting winding number are illustrated
schematically in Figure~\ref{fig:finite-mass-circuit}(b).
\begin{figure}[htbp]
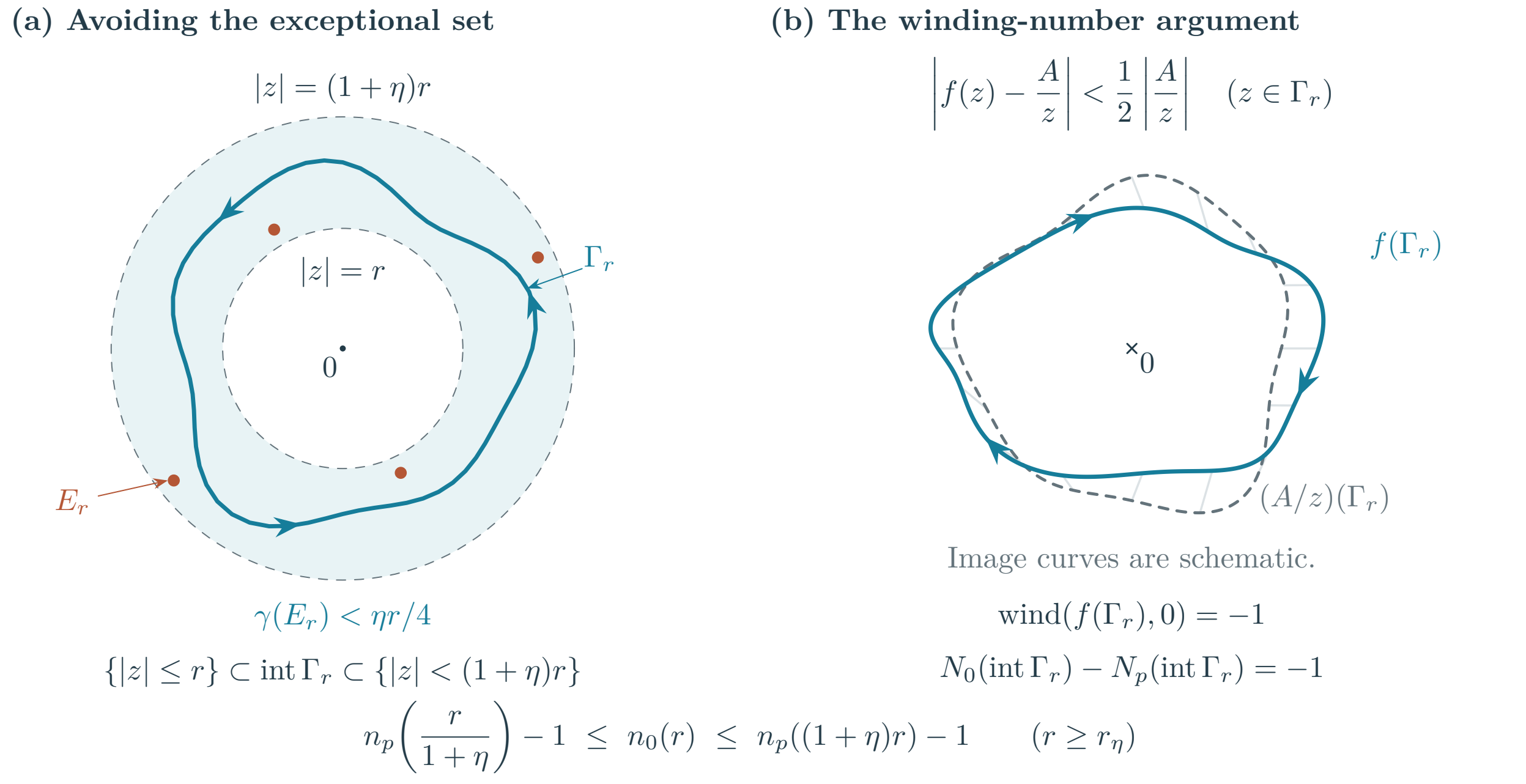

\centering
\begingroup
\resizebox{\linewidth}{!}{%
%
}
\endgroup
\caption{The circuit and winding-number argument in the proof
of Theorem~\ref{thm:finite-mass}.
(a) Small analytic capacity permits a positively oriented polygonal
Jordan curve $\Gamma_r$ in the annulus that surrounds the inner disk
and avoids the exceptional set $E_r$.
(b) Schematic image curves: the comparison with $A/z$ gives a homotopy
avoiding the origin, so the winding number of $f$ on $\Gamma_r$ is $-1$.
Here $N_0(\Omega)$ and $N_p(\Omega)$ denote the numbers of zeros and
poles in $\Omega$, counted with multiplicity. The disk inclusions then
give the counting inequalities \eqref{eq:finite-mass-counts}.}
\label{fig:finite-mass-circuit}
\end{figure}

By the argument principle,
the number of zeros of $f$
inside $\Gamma_r$ is one less
than the number of poles there,
counted with multiplicity.
The interior of $\Gamma_r$
contains the closed disk
$\{|z|\le r\}$ and is contained
in $\{|z|<(1+\eta)r\}$.
Hence
\begin{equation}\label{eq:count-two-radii}
 n_p(r)-1\le n_0((1+\eta)r),\qquad
 n_0(r)\le n_p((1+\eta)r)-1 .
\end{equation}
Replacing $r$ by $r/(1+\eta)$
in the first inequality of
\eqref{eq:count-two-radii}
gives \eqref{eq:finite-mass-counts}. Since $n_p(r)\to\infty$,
the lower bound implies that $f$ has infinitely many zeros.
This completes the proof.
\end{proof}

\end{document}